\documentclass{imsart}

\RequirePackage[OT1]{fontenc}
\RequirePackage{amssymb, amsfonts, amstext,amsthm, latexsym,amsmath}
\RequirePackage{graphicx, bm, bbm, yfonts, colortbl}
\RequirePackage[utf8]{inputenc}
\RequirePackage[numbers]{natbib}
\RequirePackage[colorlinks,citecolor=blue,urlcolor=blue]{hyperref}
\startlocaldefs
\numberwithin{equation}{section}
\numberwithin{subsection}{section}
\theoremstyle{plain}

\newtheorem{teo}{Theorem}[section]
\newtheorem{prop}{Proposition}[section]
\newtheorem{cor}{Corollary}[section]
\newtheorem{lem}{Lemma}[section]
\newtheorem{rem}{Remark}[section]
\newtheorem{exa}{Example}[section]
\newtheorem{examples}{Examples}[section]
\newcommand{\beq}{\begin{equation}}
\newcommand{\beqn}{\begin{equation*}}
\newcommand{\eeq}{\end{equation}}
\newcommand{\eeqn}{\end{equation*}}

\newcommand{\R}{\mathbb{R}}
\newcommand{\Z}{{\mathbb Z}}

\newcommand{\Var}{\mathrm{Var}}
\newcommand{\Cov}{\mathrm{Cov}}

\endlocaldefs

\begin{document}

\begin{frontmatter}
\title{Anisotropic Gaussian random fields: Criteria for hitting probabilities and applications}
\runtitle{Hitting probabilities}

\begin{aug}
\author{\fnms{Adri\'an} \snm{Hinojosa-Calleja}\thanksref{t2}\ead[label=e1]{hinojosa@ub.edu}}
\and
\author{\fnms{Marta} \snm{Sanz-Sol\'e}\thanksref{t2}\ead[label=e2]{marta.sanz@ub.edu}}

\address{Facultat de Matem\`atiques i Inform\`atica, Universitat de Barcelona\\
Gran Via de les Corts Catalanes, 585, E-08007 Barcelona, Spain\\
Barcelona Graduate School of Mathematics\\
\printead{e1,e2}}



\thankstext{t2}{Supported by the grant MTM 2015-65092-P from the \textit{Direcci\'on General de
Investigaci\'on, Ministerio de Ciencia, Innovaci\'on y Universidades}, Spain.}
\runauthor{A. Hinojosa-Calleja, M. Sanz-Sol\'e}

\affiliation{University of Barcelona}

\end{aug}

\begin{abstract}
We develop criteria for hitting probabilities of anisotropic Gaussian random fields with associated canonical pseudo-metric given by a class of gauge functions. This yields lower and upper bounds in terms of general notions of capacity and Hausdorff measure, respectively, therefore extending the classical estimates with the Bessel-Riesz capacity and the $\gamma$-dimensional Hausdorff measure. We apply the criteria to a system of linear stochastic partial differential equations driven by space-time noises that are fractional in time and either white or colored in space.
\end{abstract}

\begin{keyword}[class=MSC]
\kwd[Primary ]{60G60, 60G15, 60H15, 60H07}
\kwd{60K35}
\kwd[; secondary ]{60G22, 60G17}
\end{keyword}

\begin{keyword}
\kwd{Capacity}
\kwd{Hausdorff measure}
\kwd{hitting probabilities}
\kwd{systems of linear SPDEs}
\kwd{fractional noise}
\end{keyword}
\end{frontmatter}

\section{Introduction}
\label{s0}
This paper is a contribution to the study of hitting probabilities for anisotropic Gaussian random fields. The motivation arises from applications of probabilistic potential theory to 
systems of linear stochastic partial differential equations (SPDEs) driven by a noise fractional in time and either white or colored in space.

Let $X=\{X(x),\ x\in \R^d\}$ be a $\R^D$-valued Gaussian process with independent components. The canonical pseudo-distance corresponding to $X$ is defined by 
$\textgoth{d}(x,y) =\Vert X(x)-X(y)\Vert_{L^2(\Omega)}$.
In this article, the process $X$ is termed isotropric (repectively, anisotropic) if, up to non null multiplicative constants, $\textgoth{d}(x,y)$ is bounded below and above by  an isotropic (respectively, anisotropic) function $G$ of the variable $x-y$. We will write $\textgoth{d}(x,y)\asymp G(x-y)$. The simplest example of $G$ describing anisotropy is  
\beq
\label{s0.2}
G(x-y)=\sum_{j=1}^d |x_j-y_j|^{\alpha_j}, \quad x, y \in \R^d,\quad \alpha_j>0,
\eeq
where at least two of the $\alpha_j$'s are different. When $\alpha_j=\alpha$ for all $j$, $G$ expresses isotropy. The fractional Brownian sheet and the random field solution to  linear stochastic heat equations fall into this category of anisotropic processes, while the solution to  linear wave equations is an example of isotropic process.

The study of hitting probabilities for $X$ consists mainly in obtaining upper and lower bounds on the probabilities of random sets $F_{I,A}:=\{X^{-1}(A)\cap I\ne\emptyset\}=\{X(I)\cap A \ne\emptyset\}$, $I\subset \R^d$, $A\subset \R^D$, in terms of the Hausdorff measure and/or the capacity of the set $A$. Such estimates provide the background to characterise polarity of sets for $X$, to deduce the Hausdorff dimension of $F_{I,A}$ and, in general, to gain insight into analytic and geometric properties of the process. 

For Gaussian processes with anisotropies described by  \eqref{s0.2}, abstract results on hitting probababilities have been proved in \cite[Theorem 7.6, p. 188]{xia} and \cite[Theorem 2.1]{b-l-x}. Extensions to non Gaussian processes are proved in \cite[Theorems 2.1 and 2.4]{dal:san}. In these works, upper and lower bounds are given in terms of the usual $\gamma$-dimensional Hausdorff measure and the $\beta$-Bessel-Riesz capacity, respectively. 
There are several papers applying these results, or making use of similar approaches, to random fields that are solutions to SPDEs, either Gaussian or non Gaussian. The bounds are sharp when $\gamma=\beta$.  We refer  to \cite{dalang2018} for a representative selection of references and for a survey on the state of the art.

Stochastic heat equations driven by fractional noises  provide illustrations of anisotropic random fields with  associated pseudo-metrics not fitting the above description. Among the many examples of such equations, we will focus here in the linear SPDE studied in \cite{bal:tud} namely,
\beq
\label{s0.3}
\frac{\partial v}{\partial t} = \Delta v+\dot{W}^{H,\alpha},\quad  (t,x)\in (0 ,T]\times\mathbb{R}^d;\quad v(0,x)=v_0(x), 
\eeq
where $(\dot{W}^{H,\alpha})$ is a noise fractional in time (with Hurst parameter $H\in(1/2,1)$) and either white or colored in space (depending on the values of the parameter $\alpha$).
Details on the setting are given in Section \ref{s3} (see \eqref{s3.spde}). 

When the initial condition $v_0$ vanishes and the constraint $4H-(d-\alpha) = 2$ holds, Tudor and Xiao prove that, for any $t\in(0,T]$, 
$\Vert v(t,x)-v(t,y)\Vert_{L^2(\Omega)}\asymp\left(\log(1/|x-y|)\right)^{1/2}|x-y|$ (\cite[Theorem 4]{xia}). This result suggests the use of more general notions of Hausdorff measures and capacities than the classical $\gamma$-dimensional Hausdorff measure and the $\beta$-Bessel-Riesz capacity respectively,  to achieve sharp upper and lower bounds on the hitting probabilities for $v$.

Building on this fundamental idea, we develop our work into two steps. In the first one, we consider a class of increasing continuous real-valued functions $q$ such that $q(0)=0$, and Gaussian random fields $X$ with canonical pseudo-metric satisfying $\textgoth{d}(x,y)\asymp q(|x-y|)$. We term this type of processes {\em $q$-anisotropic}. In this setting, we establish abstract criteria for hitting probabilities. If $q$ equals the function $G$ in \eqref{s0.2}, we recover the results from \cite{dal:san}, \cite{b-l-x} and \cite{xia} mentioned above. 
The second step consists of an application of the new criteria to a system of random field solutions to \eqref{s0.3}.

After these introductory paragraphs, we describe with some detail the sections of the paper. In Section \ref{s1}, we summarize the basic notions and notations used throughout the article. Section \ref{s2} is devoted to the discussion of abstract criteria on hitting probabilities for $q$-anisotropic Gaussian processes. First, we consider the case where anisotropy of the process is described by a single function $q$ and then, we extend the analysis to the case where two (or more) different functions like $q$ are needed in the description. In the first event, we denote the processs by $M$ and call it  {\em single $q$-anisotropic} while in the second, the process is denoted by $U$ and is called {\em multiple $q$-anisotropic}. Although the first case could be deduced from the second one, for didactic reasons, we decided to take this path. However, in this introductory description, we will restrict to {\em single $q$-anisotropic} processes $M$.

The criteria for the upper bounds are proved using the strategy of \cite[Section 2]{dal:san}. The main ingredient is Lemma \ref{l2.2}, which takes the role of Lemma 2.5 in \cite{dal:san}. We remark that if the centred process $M-E(M)$ satisfies $\textgoth{d}(x,y)\le C q(|x-y|)$ then the assumptions of Lemma \ref{l2.2} hold. Applying a re-scaling defined by means of $q$, with an approach close to the proof of \cite[Theorem 2.6]{dal:san}, we deduce upper bounds for hitting probabilities for small balls (see Theorem \ref{t2.3} and Lemma \ref{l2.3}). In particular, Lemma 
 \ref{l2.3} reveals that  in our context,  the $g_q$-Hausdorff measure $\mathcal{H}_{g_q}(A)$ (see the definition in Section \ref{s1}) with $g_q(\tau)= \tau^D/(q^{-1}(\tau))$  is the suitable choice  of geometric measure for upper bounds of the hitting probabilities. The classical covering argument yields Theorem \ref{2.2-tconclusive} (see also Theorem \ref{t2.2-tconclusive-st} for the {\em multiple $q$-anisotropic} case).
We note by passing  that our results hold for processes with continuous mean function, therefore removing the constraint of being centred  in previous works.  

Recall the definition of $\textgoth{g}$-capacity given in Section \ref{s1} below. Assume that the process $M$ satisfies $\textgoth{d}(x,y) \asymp q(|x-y|)$ ($q$ is not necessarily the same function as in the preceding paragraph). 
In coherence with the classical anisotropic case, we expect the lower bounds on hitting probabilities to be given in terms of  the $(g_{ q})^{-1}$-capacity. In fact, for $\beta>0$, the $\beta$-Bessel-Riesz capacity is defined by the kernel $\tau^{-\beta}$. We prove that this is indeed the case if we restrict the class of functions $q$ for which 
$g_q(\tau)= \tau^D/(q^{-1}(\tau))$ satisfies a {\em rate growth} control at $\tau=0$. More precisely, let $\tau\mapsto v_q(\tau)$
be the radial integral of the function $(q^D(|z|))^{-1}$ over the circular ring determined by the radii $\tau$ and a constant $c_0>\tau$. We require
\beq
\label{vq-gq}
[g_q(\tau )]^{-1}= O(v_q(\tau)), \quad \tau\downarrow 0  
\eeq
(see \eqref{2.1.9}). When $q(\tau)=\tau^{\nu}$ (the classical case) this imposes no restriction (see Section \ref{appendix} for details).
Then, adding to condition \eqref{vq-gq} the set of assumptions as in the classical anisotropic case (see Hypotheses $(H_M)$ in Section \ref{s2.1}) we establish in Theorem \ref{s2.1-t1} 
criteria for lower bounds for hitting probabilities in terms of $\text{Cap}_{(g_{ q})^{-1}}(A)$. The proof combines the approach of  \cite[Theorem 2.1]{b-l-x}, \cite[Theorem 7.6]{xia}, based on weak approximations of measures, and the results of \cite[Section 3]{dal:san}. 

In the case where the function $q$ in the discussions on upper and lower bounds are equal we easily obtain that points are polar for $M$ if and only if $\lim_{\tau\downarrow 0}g_q(\tau) = 0$.

We close Section \ref{s2} with a sample of generic and concrete examples where the above criteria apply and we recover known results on the linear stochastic heat, wave and Poisson equations.

In Section \ref{s3}, we prove results on the random field solution to \eqref{s0.3} that are required in the application of the abstract criteria of Section \ref{s2}. These are on the covariance structure of the process and the identification of the associated anisotropic canonical pseudo-metric (see Lemmas \ref{l3.1} and \ref{l3.2}, and Theorem \ref{t3.1}, respectively).
This leads eventually to Theorem \ref{t3.2-ub} on sharp hitting probabilities for a system of SPDEs derived from \eqref{s0.3}. Finally, Section \ref{appendix} gathers technical details on examples where the results can be applied.

Before this work was completed, we came across the arXiv document \cite{enualart-viens}. Both articles share the aim of establishing abstract criteria on anisotropic Gaussian processes beyond the classical case. Our setting is more general and  there are many differences in the approaches. There is however similarity in 
the  proof of the criterion for the lower bound. As was mentioned above, it relies on \cite[Theorem 2.1]{b-l-x}, \cite[Theorem 7.6]{xia}.

There are several natural questions that are work in progress or in our research plans for the future. For example, the extension of the abstract criteria to non Gaussian $q$-anisotropic processes and the investigation of  applications that could provide a strong motivation for this. For instance, SPDEs driven by multiplicative fractional type noises. Deepen in the understanding of polarity is also a challenging project. In particular, what means {\em critical dimension} in the setting of this article; and then, in connection with \cite{x-m-d}, how could polarity of points at critical dimension be characterized.

\section{Preliminaries and notations}
\label{s1}

Let $g: \R_+\rightarrow \R_+$ be monotone increasing and right-continuous. Assume that on a small non empty interval $[0,\varepsilon_0]$, $g$ is strictly increasing. The {\em  $g$-Hausdorff measure} of 
 a Borel set $A\subset\mathbb{R}^D$ is defined  by 
\beqn
\mathcal{H}_g(A)=\lim \inf_{\varepsilon\downarrow 0}\left\{\sum_{i=1}^\infty g(2r_i): A\subset\bigcup_{i=1}^\infty B_{r_i}(x_i),\ \sup_{i\geq 1}r_i\leq\varepsilon\right\}
\eeqn
(see e.g. \cite{rogers}). In the particular case $g(\tau)= \tau^\gamma$, with $\gamma>0$, $\mathcal{H}_g(A)$ is the $\gamma$-dimensional Hausdorff measure, usually denoted by $\mathcal{H}_\gamma(A)$ (see e.g. \cite{mattila}).
\smallskip


A function $\textgoth{g}: \R^D\longrightarrow \R_+\cup \{\infty\}$ is a  {\em symmetric potential kernel} if: (a) $\textgoth{g}$ is symmetric; (b) $\textgoth{g}(z)> 0$, for all $z\ne 0$;  (c) $\textgoth{g}(0) = \infty$; (d) $\textgoth{g}$ is continuous on $\R^D\setminus\{0\}$. 

The energy of a measure $\mu$ on $\R^D$ relative to $\textgoth{g}$ is given by the expression
\beqn
\mathcal{\mathcal{E}}_{\textgoth{g}}(\mu) = \int_{\R^D\times \R^D} \textgoth{g}(y-\bar y)\ \mu(dy)\mu(d\bar y).
\eeqn
The $\textgoth{g}$-{\em capacity} of a Borel set $A\subset\mathbb{R}^D$ is defined  by
\beq
\label{1.1}
\text{Cap}_{\textgoth{g}}(A) = \left[ \inf_{\mu\in\mathbb{P}(A)}\mathcal{E}_{\textgoth{g}}(\mu)\right]^{-1},
\eeq
where $\mathbb{P}(A)$ denotes the set of probability measures on $A$. Since $\textgoth{g}$ is symmetric, this defines a Choquet capacity (see e.g. \cite[Theorem 2.1.1, p. 533]{kho}).

When $\textgoth{g}$ is the Bessel-Riesz kernel of order $\gamma\in \R$, the $\textgoth{g}$-{\em capacity} is  the 
 Bessel-Riesz capacity usually denoted by $\text{Cap}_\gamma (A)$ (see e.g. \cite[p. 376]{kho}).

Throughout this article, a {\em gauge function} means a strictly increasing continuous function 
$q: [0,r)\subset \R_+ \mapsto \R_+$ satisfying $q(0)=0$.

Whenever we consider the expression $\log\frac{c}{\tau}$, $\tau>0$, we assume that $c$ is large enough to ensure  $\log\frac{c}{\tau}>1$.




Throughout the paper we will use the following notations. The Euclidean norm on $\R^n$ is denoted by 
$\vert\cdot\vert$. For $x\in\R^n$ and $r\ge 0$, $B_r(x)$ denotes the open Euclidean ball centred at $x$ with radius $r$. Given $f: \R^n \rightarrow \R$, its Fourier transform is defined by the formula $\mathcal{F}f(\xi) = \int_{{\R}^n} e^{i x\cdot \xi}f(x)\ dx$, with $x\cdot \xi$ denoting the scalar product. Let $F$ be a set in a metric space $(S,d)$.  For $\rho>0$, $F^{(\rho)}$ denotes the set of points such that $d(x,F)<\rho$. Positive real constants are generically denoted by the letter $C$, or variants, like $\bar C$, $\tilde C$, $c$, etc. If we want to make explicit the dependence on some parameters  $a_1, a_2, \ldots$,
we write $C(a_1, a_2, \ldots)$ or  $C_{a_1, a_2, \ldots}$. The symbol
$\asymp$ between two mathematical expressions means equivalence up to multiplicative constants.

\section{Criteria for hitting probabilities}
\label{s2} 
We devote this section to investigate hitting probabilities of Gaussian random fields. The main results are Theorems \ref{2.2-tconclusive}, \ref{t2.2-tconclusive-st}, \ref{s2.1-t1} and \ref{s2.1-t2}, which yield upper and lower bounds in terms of the notions of $g$-Hausdorff measure and $\textgoth{g}$-capacity, respectively. 

\subsection{Upper bounds for hitting probabilities}
\label{s2.2}
The aim of this subsection is to prove extensions of Theorem 2.6 in \cite{dal:san} on sufficient conditions for upper bound estimates of hitting probabilities of Gaussian processes. 
\medskip

\noindent{\bf The single $q$-anisotropic case}
\medskip

The results of the first part of the section concern a $D$-dimensional stochastic process denoted by
\beq
\label{defM}
M=\{M(x)= (M_1(x), \ldots, M_D(x)),\  x\in \R^d\}.
\eeq
We start with a technical lemma which is a generalised version of \cite[Lemma 2.5]{dal:san}.

\begin{lem}
\label{l2.2} 
Assume that the process $M$ has continuous sample paths a.s. Let $q:\R_+\longrightarrow \mathbb{R}_+$ be a differentiable gauge function. Suppose that for all $\varepsilon\in (0,1)$ small enough and $x\in\R^d$,
\beq
\label{eq2.2.17}
E\left(\int_{B_\varepsilon(x)}dy \int_{B_\varepsilon(x)}d\bar y\  \exp\left(\frac{\vert M(y) -M(\bar y) \vert}{q(\vert y-\bar y\vert)}\right) \right)\leq C\varepsilon^{2d},
 \eeq
 for some constant $C$. Set $S_\varepsilon(x) = B_{\frac{q^{-1}(\varepsilon)}{2}}(x)$. Then, the following statements hold.
 \begin{enumerate}
\item For all $p\geq 1$, there exist constants $C(p,d)$ and  $\tilde C(d)$ such that for $\varepsilon$ small enough,
\begin{align}
\label{eq2.2.18}
&E\bigg(\sup_{y\in S_\varepsilon(x)} \vert M(y)-M(x) \vert ^p \bigg)\notag\\
&\qquad\quad \leq  
C(p,d)\  \varepsilon^{p-1}  q^{-1}(\varepsilon)
 \int_0^1\log^p\left(1+\frac{\tilde C(d)}{\tau^{2d}}\right) \dot q\left(q^{-1}(\varepsilon)\tau\right)\ d\tau,
\end{align} 
where $\dot q$ denotes the derivative of $q$.
\item Assume that $q$ is such that, for any $r, \tau\in[0,c_0]$, with $c_0>0$ sufficiently small,
\beq
\label{hq}
q(r\tau) \le  \varphi(\tau) q(r), \ \ \dot q(r\tau) \le \frac{1}{r} \psi (\tau)q(r\tau),
\eeq
where $\varphi$ and $\psi$ are Borel functions such that, denoting $\Phi(\tau)=\varphi(\tau)  \psi (\tau)$, we have
\beq
\label{logfinite}
\int_0^1\log^p\left(1+\frac{\tilde C(d)}{\tau^{2d}}\right) \Phi(\tau)\ d\tau < \infty.
\eeq
Then, for all $p\geq 1$, there exists a constant $C(p,d)$ such that for all $\varepsilon$ small enough,
\beq
\label{eq2.2.19}
E\bigg(\sup_{y\in S_\varepsilon(x)} \vert M(y)-M(x) \vert ^p \bigg)\leq C(p,d)\varepsilon^p.
\eeq
\end{enumerate}
\end{lem} 
\proof 
1. Let
\beq
\label{cepsilon}
\mathcal{C}_\varepsilon(\omega)=\int_{S_{\varepsilon}(x)}dy \int_{S_{\varepsilon}(x)}d\bar y\  \exp\left(\frac{\vert M(y,\omega)-M(\bar y,\omega)\vert}{q(\vert y-\bar y\vert) }\right).
\eeq 
From \eqref{eq2.2.17}, we deduce $\mathcal{C}_\varepsilon (\omega) < \infty$, a.s. Notice that for almost all $\omega$, $\mathcal{C}_\varepsilon (\omega) \ge C_2 \left(q^{-1}(\varepsilon)\right)^{2d}$, for some constant $C_2>0$.

Applying \cite[Proposition A.1, (A.3)]{dal:kho} to $S:=S_{\varepsilon}(x)$ endowed with the Euclidean distance $\rho$, $\mu$ there the Lebesgue measure, $\Psi(\tau): = e^\tau-1$ and $p(\tau):=q(\tau)$, we deduce
\beqn
\sup_{y\in S_{\varepsilon}(x) }\vert M(y)-M(x)\vert 
\le 10 \int_0^{q^{-1}(\varepsilon)} \Psi^{-1}\left(\frac{C_1 \mathcal{C}_\varepsilon(\omega)}{\tau^{2d}}\right) \dot q(\tau)\ d\tau,
\eeqn
with $C_1$ depending on $d$. Here, we have used that the volume of the $d$-dimensional Euclidean ball of radius $r$ equals a multiple constant times $r^d$. Therefore, for any $p\ge 1$,
\begin{align}
\label{2.2.17bis}
&E\left(\sup_{y\in S_{\varepsilon}(x) }\vert M(y)-M(x)\vert^p \right)
\le 10^pE\left(\left\vert \int_0^{q^{-1}(\varepsilon)} \Psi^{-1}\left(\frac{C_1 \mathcal{C}_\varepsilon(\omega)}{\tau^{2d}}\right) \dot q(\tau)\ d\tau\right\vert^p\right)\notag\\
&\qquad\qquad\le 10^p \left(q\left(q^{-1}(\varepsilon)\right)\right)^{p-1}E\left(\int_0^{q^{-1}(\varepsilon)} 
\log^p\left(1+\frac{C_1 \mathcal{C}_\varepsilon(\omega)}{\tau^{2d}}\right) \dot q(\tau)\ d\tau\right)\notag\\
&\qquad\qquad= C(p) \varepsilon^{p-1}\int_0^{q^{-1}(\varepsilon)} 
E\left[\log^p\left(1+\frac{C_1 \mathcal{C}_\varepsilon(\omega)}{\tau^{2d}}\right)\right] \dot q(\tau)\ d\tau,
\end{align}
where in the second inequality, we have applied H\"older's inequality with respect to the measure $\dot q(\tau) d\tau$. Observe that we may take $C_1$ as large as we want.

The function $x\mapsto \log^p(1+x)$ is concave on $[e^{p-1}-1,\infty)$. Hence, by taking $C_1\ge (e^{p-1}-1) C_2^{-1}$, we can apply Jensen's inequality to estimate from above the term 
$E\left[\log^p\left(1+\frac{C_1 \mathcal{C}_\varepsilon(\omega)}{\tau^{2d}}\right)\right]$ on the right-hand side of \eqref{2.2.17bis}. By doing so, then
using  \eqref{cepsilon} and \eqref{eq2.2.17}, and applying the change of variables $\tau\mapsto \left(q^{-1}(\varepsilon)\right)^{-1}\tau$, we obtain, 
\begin{align}
\label{2.2.17tris}
&E\left(\sup_{y\in S_{\varepsilon}(x) }\vert M(y)-M(x)\vert^p \right)\notag\\
&\qquad \qquad\le C(p,d)\  \varepsilon^{p-1}  q^{-1}(\varepsilon)
 \int_0^1\log^p\left(1+\frac{\tilde C(d)}{\tau^{2d}}\right) \dot q\left(q^{-1}(\varepsilon)\tau\right)\ d\tau,
\end{align}
with some constant $\tilde C$ depending on $d$. This ends the proof of \eqref{eq2.2.18}.
\medskip

2. The conditions \eqref{hq} imply $r\dot q(r\tau)\le \Phi(\tau)q(r)$. For $r:=q^{-1}(\varepsilon)$ this yields
\beqn
q^{-1}(\varepsilon)\dot q(q^{-1}(\varepsilon) \tau)\le \Phi(\tau)\varepsilon.
\eeqn
Thus, up to the multiplicative constant $C(p,d)$, the right-hand side of \eqref{eq2.2.18} is  equal to \break
$\varepsilon^p \int_0^1\log^p\left(1+\frac{\tilde C(d)}{\tau^{2d}}\right) \Phi(\tau)\ d\tau$ and therefore, assuming \eqref{logfinite}, we obtain \eqref{eq2.2.19}.
\qed

\begin{examples}
\label{2.2-e0}
We exhibit two examples of gauge functions $q$ that satisfy the hypotheses of Lemma \ref{l2.2}.
\begin{enumerate}
\item $q(\tau) = \tau^\nu$, $\tau>0$, $\nu>0$.  The conditions \eqref{hq} hold for any $\tau, r >0$, with
$ \varphi(\tau) = \tau^\nu, \  \psi(\tau) = \frac{\nu}{\tau}$.
Since $\int_0^1\log^p\left(1+\frac{\tilde C(d)}{\tau^{2d}}\right) \tau^{\nu-1}\ d\tau<\infty$ for any $p\ge 1$, condition \eqref{logfinite} holds.
 \item $q(\tau) = \tau^\gamma\left(\log\frac{c}{\tau}\right)^\delta$, $\tau>0$, with $\gamma>0$, $\delta\ge 0$. Then,
 \begin{align*}
 q(r\tau) &= r^\gamma\tau^\gamma\left(\log\frac{c}{r\tau}\right)^\delta
 \le r^\gamma\tau^\gamma\left(\log\frac{C}{r}+ \log\frac{C}{\tau}\right)^\delta\\
 &\le C(\delta)r^\gamma\tau^\gamma\left(\log\frac{C}{r}\right)^\delta \left(1+\log\frac{C}{\tau}\right)^\delta,
 \end{align*}
 with $C^2\ge c$. Hence,
 \beq
 \label{exq2-1}
  q(r\tau)\le \varphi(\tau)q(r), \ {\text{with}}\ \   \varphi(\tau) =  C(\delta)\tau^\gamma\left(1+\log\frac{C}{\tau}\right)^\delta.
  \eeq
The derivative of $q$ is 
$\dot{q}(\tau)=\tau^{\gamma-1}\left(\log \frac{c}{\tau}\right)^{\delta-1} (\gamma\log \frac{c}{\tau}-\delta)$ and therefore,
it is increasing on $[0,c e^{-\frac{\delta}{\gamma}}]$. In the sequel, we will restrict $q$ to this interval, therefore
$
\dot{q}(\tau)\le \gamma \tau^{\gamma-1}\left(\log \frac{c}{\tau}\right)^{\delta} \le \gamma\frac{1}{\tau}q(\tau).
$
Consequently, 
  \beq
  \label{exq2-2}
  \dot q(r\tau)\le \frac{1}{r} \psi(\tau)q(r\tau), \  {\text{with}}\  \psi(\tau) = \frac{\gamma}{\tau}.
  \eeq
Since $\Phi(\tau) = C(\delta,\gamma)\tau^{\gamma-1}\left(1+\log\frac{C}{\tau}\right)^\delta$, we see that condition \eqref{logfinite} holds.
\end{enumerate}
\end{examples}

\begin{rem}
\label{s2.2-r10}
Let $q$ be a function as in Lemma \ref{l2.2}. Assume that the process $M$ in Lemma \ref{l2.2} is Gaussian, centred and such that, there exists a constant $C$ and for any $|y-\bar y| < 2\varepsilon$,
\begin{equation}
\label{1952}
\Vert M(y) - M(\bar y)\Vert_{L^2(\Omega)} \le C q(|y - \bar y|).
\end{equation}

Then $M$ satisfies the condition \eqref{eq2.2.17} for any $x\in\R^d$. Indeed, \eqref{1952} implies
\begin{align*}
\exp\left(\frac{|M(y) -M(\bar y)|}{q(|y - \bar y|)}\right) &\le \exp\left(\frac{1}{C} \frac{|M(y) -M(\bar y)|}{\sqrt{{\text{Var}}\ (M(y) -M(\bar y))}}\right) = \exp(c |Z|),
\end{align*}
where $c=1/C$ and $Z$ is a standard Gaussian random variable. Since, $E([exp(c |Z|)])$ is finite, \eqref{eq2.2.17} holds.
\end{rem}
\medskip

For any $\varepsilon\in(0,1)$, $j\in\mathbb{Z}^d$, $j=(j_1,\ldots,j_d)$,  set
\beq
\label{r}
R^{\varepsilon}_j=\prod_{i=1}^{d}\left[\frac{q^{-1}(\varepsilon)}{\sqrt{d}}j_i, \frac{q^{-1}(\varepsilon)}{\sqrt{d}}(j_i+1)\right],
\eeq
and for $x\in R^{\varepsilon}_j$, define $ x^{\varepsilon}_j:=\left(\frac{q^{-1}(\varepsilon)}{\sqrt{d}}j_i\right)_{i=1,...,d}$. 
Observe that ${\text{diam}}(R^{\varepsilon}_j) = q^{-1}(\varepsilon)$ and $R^{\varepsilon}_j\subset B_{\frac{q^{-1}(\varepsilon)}{2}}(\bar x^{\varepsilon}_j)$, where $\bar x^{\varepsilon}_j= \left(\frac{q^{-1}(\varepsilon)}{\sqrt{d}}(j_i+\frac{1}{2})\right)_{i=1,...,d}$. Moreover, by the triangle inequality,
\beq
\label{randb}
\sup_{x\in R^{\varepsilon}_j}\left(\left\vert M(x)-M(x_j^\varepsilon)\right\vert\right)\le 2\sup_{x\in B_{\frac{q^{-1}(\varepsilon)}{2}}(\bar x^{\varepsilon}_j)}\left(\left\vert M(x)-M(\bar x^{\varepsilon}_j)\right\vert\right).
\eeq
\smallskip

The next statement provides an extension of \cite[Theorem 2.6, (26)]{dal:san} to non necessarily centred processes. 
\begin{teo}\label{t2.3}
Fix $\varepsilon\in(0,1)$ small enough, $j=(j_1,\ldots,j_d)\in\mathbb{Z}^d$, and 
let $R_j^\varepsilon$ be as in \eqref{r}. Assume that the process $M$ given in \eqref{defM} is Gaussian, continuous, with i.i.d. components and such that  $\sigma^2:=\Var(M_1(x_j^\varepsilon))>0$. 

1. Let $f(x) = E(M(x))$ and $\tilde M(x)=M(x)- f(x)$. We assume that  $f:\R^d\rightarrow \R^D$ is a continuous map and furthermore, for some constant $C(d,D)$,
 \begin{align}
\label{eq2.2.7}
E\left( \sup_{x\in R^{\varepsilon}_{j}}\left\vert \tilde M(x)-\tilde M(x^{\varepsilon}_j)\right\vert^2\right)&\leq C(d,D)\varepsilon^2.
\end{align}
 Then there exists a constant $C(\sigma,d,D)$  such that, for every $z\in\mathbb{R}^D$,
\begin{equation}
\label{eq2.2.8}
P\left(M(R^{\varepsilon}_{j})\cap B_\varepsilon(z)\neq \emptyset\right)\leq C(\sigma,d,D)\varepsilon^D .
\end{equation}

2. Suppose that for some constant $\bar C(d,D)$
 \begin{align}
\label{eq2.2.7-bis-bis}
E\left( \sup_{x\in R^{\varepsilon}_{j}}\left\vert M(x)-M(x^{\varepsilon}_j)\right\vert^2\right)&\leq \bar C(d,D)\varepsilon^2.
\end{align}
Then there exists a constant $\bar C(\sigma,d,D)$  such that, for every $z\in\mathbb{R}^D$,
\begin{equation}
\label{eq2.2.8-bis}
P\left(M(R^{\varepsilon}_{j})\cap B_\varepsilon(z)\neq \emptyset\right)\leq \bar C(\sigma,d,D)\varepsilon^D .
\end{equation}
\end{teo}
\proof 
 1. We follow the approach of \cite[Theorem 2.6]{dal:san} with some modifications due to the fact that the process $M$ is not centred.


Because $M$ is continuous, for any $z\in\mathbb{R}^D$ we have
\begin{equation*}
P\left(M(R^{\varepsilon}_{j})\cap B_{\varepsilon}(z)\neq \emptyset\right)=P\left(\inf_{x\in R^{\varepsilon}_{j}}\vert M(x)-z\vert\leq \varepsilon\right).
\end{equation*}

Assume we can prove that there exists a constant $c(\sigma,d,D)$ such that for any $z_1\in\mathbb{R}$,
\beq
\label{one-dimensional}
P\left(\inf_{x\in R^{\varepsilon}_{j}}\vert M_1(x)-z_1\vert\leq \varepsilon\right) \le c(\sigma,d,D)\ \varepsilon,
\eeq
where $M_1(x)$ is the first component of the random vector $M(x)$. Then, because the components of $M(x)$ are i.i.d, \eqref{one-dimensional} yields 
\eqref{eq2.2.8} with $C(\sigma,d,D) = [c(\sigma,d,D)]^D$.

For the proof of \eqref{one-dimensional}, we fix $x\in R^{\varepsilon}_{j}$ and compute the conditional expectation
\begin{equation}
\label{eq2.2.11}
E\left(M_1(x)\vert \tilde{M}_1(x^{\varepsilon}_j)\right)= f_1(x)+E\left(\tilde M_1(x)\vert \tilde M_1(x^{\varepsilon}_j)\right)
=f_1(x)+ c^{\varepsilon}_{j}(x)\tilde M_1(x^{\varepsilon}_j),
\end{equation}
where
\beqn
c^{\varepsilon}_{j}(x)=\frac{\Cov\left(\tilde{M}_1(x), \tilde{M}_1(x^{\varepsilon}_j)\right)}{\Var \left(\tilde{M}_1(x^{\varepsilon}_j)\right)}.
\eeqn 

Define
\beqn
Y^{\varepsilon}_{j}=\inf_{x\in R^{\varepsilon}_{j}}\left\vert E\left({M}_1(x)\vert \tilde M_1(x^{\varepsilon}_j)\right)-z_1\right\vert,\quad
Z^{\varepsilon}_{j}=\sup_{x\in R^{\varepsilon}_{j}}\left\vert {M_1}(x)- E\left({M}_1(x)\vert \tilde M_1(x^{\varepsilon}_j)\right)\right\vert.
\eeqn 
These are independent random variables satisfying
\begin{equation}
\label{eq2.2.10}
P\left( \inf_{x\in R^{\varepsilon}_j}\vert {M}_1(x)- z_1\vert\leq \varepsilon \right)\leq P\left(Y^{\varepsilon}_{j}\leq \varepsilon+Z^{\varepsilon}_{j}\right).
\end{equation}

We next prove that, for any $ r\geq 0$, 
\begin{equation}
\label{eq2.2.14}
 P(Y^{\varepsilon}_{j}\leq r)\leq C(\sigma,d,D) r.
\end{equation}
As an auxiliary result for this, we first check that for  all $\varepsilon>0$ small enough and $x\in R^{\varepsilon}_j$,
 \beq
 \label{eq2.2.12}
\vert c^{\varepsilon}_{j}(x)-1\vert \leq C(\sigma,d,D) \varepsilon,
\eeq
implying that, 
for all $\varepsilon>0$ small enough, say $\varepsilon\le \varepsilon_0(\sigma,d,D)$, and for all $x\in R^{\varepsilon}_j$, 
we have
\beq
\label{eq2.2.13}
c^{\varepsilon}_j(x) \geq \frac{1}{2}.
\eeq
 Indeed,  because $\Var(\tilde{M}_1(x_j^\varepsilon))=\Var({M}_1(x_j^\varepsilon))=\sigma^2>0$, using (\ref{eq2.2.7}), similarly as in \cite[(30), p.1356]{dal:san}, we deduce
\begin{align*}
\left\vert c^{\varepsilon}_{j}(x)-1\right\vert
\leq  \left(\frac{E\left[\tilde{M}_1(x^{\varepsilon}_j)-\tilde{M}_1(x)\right]^2}{\Var (\tilde{M}_1(x^{\varepsilon}_j))}\right)^{\frac{1}{2}}
\leq C(\sigma,d,D)\varepsilon.
\end{align*}


By \eqref{eq2.2.11}, 
\beqn
\{Y^{\varepsilon}_{j}\leq r)\} = \{\inf_{x\in R^{\varepsilon}_{j}}\vert f_1(x)+ E(\tilde{M}_1(x)\vert \tilde{M}_1(x^{\varepsilon}_j))-z_1\vert\leq r\} 
\eeqn
and the inequality $\left\vert f_1(x)+E\left(\tilde{M}_1(x)\vert \tilde{M}_1(x^{\varepsilon}_j)\right)-z_1\right\vert\leq r$ is equivalent to
\beqn
\frac{z_1-f_1(x)}{c^{\varepsilon}_{j}(x)}-\frac{r}{c^{\varepsilon}_{j}(x)}\leq \tilde{M}_1(x^{\varepsilon}_j)\leq\frac{z_1-f_1(x)}{c^{\varepsilon}_{j}(x)}+\frac{r}{c^{\varepsilon}_{j}(x)}.
\eeqn
Since by \eqref{eq2.2.13},   $\inf_{x\in R_{j}^{\varepsilon}} c^{\varepsilon}_{j}(x)\ge \frac{1}{2} $, the above remarks yield
\begin{align}
\label{eq2.2.16}
P(Y^{\varepsilon}_{j}\leq r) &\leq \sup_{s\in\mathbb{R}}
P\left( s-2r\le  \tilde{M}_1(x^{\varepsilon}_j)\leq s+ 2r\right)
=\sup_{s\in\mathbb{R} }P\left( \tilde{M}_1(x^{\varepsilon}_j)\in B_{2r}(s) \right).
\end{align}
Because the density of $\tilde{M}_1(x^{\varepsilon}_j)$ is bounded by $(\Var(M_1(x_j^\varepsilon))2\pi)^{-1/2}=\tfrac{1}{\sigma \sqrt{2\pi}}$, we have
\beqn
P(Y^{\varepsilon}_{j}\leq r) \le \sup_{s\in\mathbb{R}}P\left( \tilde{M}_1(x^{\varepsilon}_j)\in B_{2r}(s) \right)
\le C(\sigma)\ r.
\eeqn
This proves \eqref{eq2.2.14}.
\smallskip

We now address the last step in the proof of \eqref{one-dimensional}. 
From \eqref{eq2.2.10}, and because $Y^{\varepsilon}_{j}$ and $Z^{\varepsilon}_{j}$ are independent, by using \eqref{eq2.2.14}
we obtain,
\begin{align}
\label{newnew}
P\left( \inf_{x\in R^{\varepsilon}_j}\vert {M}_1(x)- z\vert\leq\varepsilon\right)&\leq c(\sigma,d,D)\ E\left[\left(\varepsilon+Z^{\varepsilon}_{j}\right)\right]\notag\\
& = c(\sigma,d,D)\left[\varepsilon+E\left(\left(Z^{\varepsilon}_{j}\right)\right)\right].
\end{align}
Since
${M}_1(x)-E\left({M}_1(x)\vert \tilde{M}_1(x^{\varepsilon}_j)\right)=\tilde{M}_1(x)-c^{\varepsilon}_{j}(x)\tilde{M}_1(x^{\varepsilon}_j)$ (see \eqref{eq2.2.11}), by the triangle inequality we have  $Z^{\varepsilon}_{j}\leq Z^{\varepsilon}_{j,1}+Z^{\varepsilon}_{j,2}$, with
\beqn
Z^{\varepsilon}_{j,1}=\sup_{x\in R^{\varepsilon}_j}\left\vert \tilde{M}_1(x)-\tilde{M}_1(x^{\varepsilon}_j)\right\vert, \quad
Z^{\varepsilon}_{j,2}=\sup_{x\in R^{\varepsilon}_j}\left\vert 1-c^{\varepsilon}_j(x)\right\vert\left\vert \tilde{M_1}(x^{\varepsilon}_j)\right\vert.
\eeqn
Apply (\ref{eq2.2.7})  to obtain  
$E(Z^{\varepsilon}_{j,1})\leq C(d,D)\varepsilon$.   
Also, as a consequence of  \eqref{eq2.2.7} and \eqref{eq2.2.12}, we have $E(\vert Z^{\varepsilon,2}_j\vert^D)\leq C(\sigma,d,D)\varepsilon$. This yields 
$E\left[\left(Z^{\varepsilon}_{j}\right)\right]\leq C(\sigma,d,D)\varepsilon$.
Along with \eqref{newnew}, this implies \eqref{one-dimensional} and, as was argued above, the proof of claim 1 is complete.

2. With minor changes in the previous proof, we can check that claim 2 holds. The details are left to the reader.
\qed
\begin{rem}
\label{match}
 In the setting of Theorem \ref{t2.3}, suppose in addition that the process $M$ is continuous (which implies that $\tilde M$ is continuous too).
Assume that $\tilde M$ satisfies the hypotheses of Lemma \ref{l2.2} with $x:=\bar x_j^\varepsilon$. Then, applying \eqref{randb} with $M$ there replaced by $\tilde M$, we see that condition \eqref{eq2.2.7}  holds for any $p\ge1$. 
 Similarly, if $M$ satisfies the hypotheses of Lemma \ref{l2.2} with $x:=\bar x_j^\varepsilon$. Then applying \eqref{randb}, we deduce that condition \eqref{eq2.2.7-bis-bis} holds for any $p\ge1$. 
\end{rem}
\smallskip

\noindent{\em Hitting probabilities for small balls}
\medskip

Using a standard argument based on total probabilities, we can derive upper bounds for hitting probabilities of small balls in terms of the function $q$, as follows.
\begin{lem}
\label{l2.3}
Let $K\subset \R^d$ be a compact set of positive Lebesgue measure. Fix $z\in\R^D$ and $\varepsilon>0$ (small enough). Let $M$ be the process defined in \eqref{defM} and assume that it is Gaussian and continuous, with i.i.d. components, and such that $\sigma^2_K:=\inf_{x\in K^{(\eta)}} {\rm{Var}}\ (M(x))>0$ (for some $\eta$ sufficiently small).
Let $q$ be a function satisfying the conditions of Lemma \ref{l2.2}, and assume that $M$ satisfies \eqref{eq2.2.17} for any $x\in K^{(\eta)}$.

Then, there exists a constant $C(K,\sigma_K,d,D)$ such that,
\beq
\label{sbg}
P\left(M(K)\cap B_\varepsilon(z)\ne \emptyset\right) \le C(K,\sigma_K,d,D) \frac{\varepsilon^D}{\left(q^{-1}(\varepsilon)\right)^d}.
\eeq
\end{lem}
\proof
Since $K$ is compact, there is a finite number of sets $R_j^\varepsilon$ (defined in \eqref{r}) satisfying $K\cap R_j^\varepsilon\neq \emptyset$; this number is a constant (depending on the dimension $d$) multiple of $\left(\frac{q^{-1}(\varepsilon)}{\sqrt{d}}\right)^{-d}$. Moreover, by Lemma \ref{l2.2} and the inequality \eqref{randb}, we see that the condition \eqref{eq2.2.7-bis-bis} holds for any $R_j^\varepsilon$ such that $K\cap R_j^\varepsilon\ne \emptyset$ and this implies \eqref{eq2.2.8-bis}. Thus,
\begin{align}
\label{sbg-bis}
&P\left(M(K)\cap B_\varepsilon(z)\ne \emptyset\right) \le \sum_{j\in\Z^d: K\cap R_j^\varepsilon\neq \emptyset}
P\left(M(R_j^\varepsilon)\cap B_\varepsilon(z)\ne \emptyset\right)\notag\\ 
&\qquad\qquad\le \tilde C(K,\sigma_K,d,D)\  \varepsilon^D \left(\frac{q^{-1}(\varepsilon)}{\sqrt{d}}\right)^{-d}
= C(K,\sigma_K,d,D)\  \varepsilon^D \left(q^{-1}(\varepsilon)\right)^{-d}.
\end{align}
\qed

For a gauge function $q$, define
\beq
\label{gg}
g_q(\tau) = \frac{\tau^D}{\left(q^{-1}(\tau)\right)^d}, \quad \tau\in\R_+.
\eeq

From Lemma \ref{l2.3} we deduce conditions for points to be polar, as follows. 

\begin{cor}
\label{s2-c10}
The hypotheses are as in Lemma \ref{l2.3}. Assume further that
\beq
\label{2.c-polar}
\lim_{\tau\downarrow 0}g_q(\tau) = 0.
\eeq
Then, for any $z\in\R^D$, $P\left(M(K)\cap \{z\}\ne \emptyset\right)= 0$, that is $\{z\}$ is polar for the process $M$ restricted to $K$.
\end{cor}
\proof
For any $\varepsilon>0$, we have
$P\left(M(K)\cap \{z\}\ne \emptyset\right)\le P\left(M(K)\cap B_\varepsilon(z)\ne \emptyset\right)$.
Applying \eqref{sbg-bis} and using \eqref{2.c-polar} yields the result.
\qed
\medskip

\noindent{\em Covering argument}
\medskip

Assume that the hypotheses of Lemma \ref{l2.3} hold. Let $g_q$ be the function defined in \eqref{gg} and assume that on a sufficiently small interval $(0,\rho_0)$, $g_q$ is strictly increasing.  
Fix $\varepsilon$ small enough. By the definition of the Hausdorff $g_q$-measure $\mathcal{H}_{g_q}(A)$, there exists a sequence of balls $(B_i, i\ge 1)$ with radii $r_i\in(0,\varepsilon)$, such that $B_i\cap A \neq \emptyset$,
$A\subset \cup_{i\ge 1}B_i$, and
$\sum_{i\ge 1} g_q(2r_i)\le \mathcal{H}_{g_q}(A)+\varepsilon$.
Then from \eqref{sbg}, for any Borel set $A\subset \R^D$ we deduce,
\begin{align*}
P\left(M(K)\cap A\ne \emptyset\right) &\le \sum_{i\ge 1} P\left(M(K)\cap B_i\ne \emptyset\right)\\
&\le C(K,\sigma_K,d,D) \sum_{i\ge 1} g_q(2r_i)
\le \mathcal{H}_{g_q}(A)+\varepsilon.
\end{align*}
Letting $\varepsilon$ tend to zero, we obtain
\beq
\label{hg}
P\left(M(K)\cap A\ne \emptyset\right)\le C(K,\sigma_K,d,D)\ \mathcal{H}_{g_q}(A).
\eeq

\bigskip

\noindent{\em Hitting probabilities in terms of $g$-Hausdorff measures}
\medskip

We summarise the preceding discussion in the following statement.
\begin{teo}
\label{2.2-tconclusive}
Let $K\subset \R^d$ be a compact set of positive Lebesgue measure. Consider a Gaussian continuous stochastic process 
$M=\{M(x)= (M_1(x),\ldots, M_D(x)), x\in  \R^d\}$ with i.i.d. components and such that $\sigma_K^2:=\inf_{x\in K^{(\eta)}} {\rm{Var}}\ (M(x))>0$
(for some $\eta>0$ sufficiently small).
Let $q$ be a function satisfying the hypotheses of Lemma \ref{l2.2} and such that the function $g_q$ given in \eqref{gg} is strictly increasing on a small interval $(0,\rho_0)$. Assume also that the process $M$ satisfies the condition \eqref{eq2.2.17} for any $x\in K^{(\eta)}$.

Then there exists a constant $C(K,\sigma_K,d,D)$ such that for any Borel set $A\subset \R^D$, 
\beq
\label{hitting-up-1}
P\left(M(K)\cap A\ne \emptyset\right)\le  C(K,\sigma_K,d,D)\ \mathcal{H}_{g_q}(A).
\eeq
\end{teo}

In Lemma \ref{s2.2-l-grow-g-st} (sections 1. and 3.)  examples of gauge functions $q$ satisfying the assumptions of Theorem \ref{2.2-tconclusive} and \eqref{2.c-polar} are given.
\newpage


\noindent{\bf The multiple $q$-anisotropic case}
\medskip

Through condition \eqref{eq2.2.17}, the function $q$ in Lemma \ref{l2.2} provides a control of the oscillations of the sample paths of the process $M$. However, as we will see in Section \ref{s3}, 
there exist stochastic processes where two or a finite number of distinct functions $q$ are needed for such control. In the second part of this section, we develop an extension of the previous results in a setting suitable for their application in Section \ref{s3}.

Let
\beq
\label{defU}
U=\{U(t,x)=(U_1(t,x),\ldots,U_D(t,x)), (t,x)\in\mathbb{R}^{d_1}\times \mathbb{R}^{d_2}\},
\eeq
be a stochastic process. In the examples in mind, the parameter $t$ refers to time (therefore $d_1=1$) while $x$ refers to space. 

\begin{lem}
\label{l2.2-ts}
Assume that the process $U$ has continuous sample paths a.s. Let $q_1$, $q_2$ be functions satisfying the properties of $q$ in Lemma \ref{l2.2}. Fix compact sets $I\subset\mathbb{R}^{d_1}$, $J\subset\mathbb{R}^{d_2}$ of positive Lebesgue measure and assume that,
for any $\varepsilon$ small enough,
\begin{align}
\label{qoneqtwo}
&E\left(\sup_{x\in J}\int_{B_\varepsilon(t)} ds \int_{B_\varepsilon(t)} d\bar s\  \exp\left(\frac{\vert U(s,x)-U(\bar s,x)\vert}{q_1(|s-\bar s|)}\right)\right) \le C \varepsilon^{2d_1},\quad t\in I, \notag\\
&E\left(\sup_{t\in I}\int_{B_\varepsilon(x)} dy \int_{B_\varepsilon(x)} d\bar y\ \exp\left(\frac{\vert U(t,y)-U(t,\bar y)\vert}{q_2(|y-\bar y|)}\right)\right) \le C \varepsilon^{2d_2},\quad x\in J.
\end{align}
Let $S_\varepsilon^1(t)=B_{\frac{q_1^{-1}(\varepsilon)}{2}}(t)$, $S_\varepsilon^2(x)=B_{\frac{q_2^{-1}(\varepsilon)}{2}}(x)$ and
$\tilde S_\varepsilon(t,x) = S_\varepsilon^1(t)\times S_\varepsilon^2(x)$. 
Then, for all $p\ge 1$, there exists a constant $C(p,d_1,d_2)$ such that, for all $\varepsilon$ small enough and $(t,x)\in I\times J$,
\beq
\label{double-st}
E\left(\sup_{(s,y)\in \tilde S_\varepsilon(t,x)}|U(s,y)-U(t,x)|^p\right) \le C(p,d_1,d_2)\ \varepsilon^p.
\eeq
\end{lem}
\proof
 We follow the steps of the proof of Lemma \ref{l2.2} considering first the processes $M:= U^{(x)}= \{U(t,x), t\in I\}$, and then 
$M:=U^{(t)}=\{U(t,x), x\in J\}$, obtained from $U$ by fixing the indices $x\in J$ and $t\in I$, respectively. The hypotheses \eqref{qoneqtwo} play 
the role of \eqref{eq2.2.17} in Lemma \ref{l2.2} for the proof of \eqref{eq2.2.19}. In this way, we obtain,
\begin{align*}
E\left(\sup_{x\in J}\sup_{s\in S_\varepsilon^1(t)} \left|U(s,x)-U(t,x)\right\vert^p\right)& \le C(p,d_1)\varepsilon^p,\\
E\left(\sup_{t\in I}\sup_{y\in S_\varepsilon^2(x)} \left|U(t,y)-U(t,x)\right\vert^p\right)& \le C(p,d_2)\varepsilon^p, 
\end{align*}
for some constants $C(p,d_1)$, $C(p,d_2)$. Using the triangle inequality, we deduce \eqref{double-st}.
\qed
\begin{rem}
\label{s2.2-r20}
Let $q_1$, $q_2$ be functions as in Lemma \ref{l2.2-ts}. Assume that the process $U$ in Lemma \ref{l2.2-ts} is Gaussian, centred and such that for any $|s-\bar s|<2\varepsilon$, $|y-\bar y|<2\varepsilon$,
\begin{align}
\label{s2.2-200}
\sup_{x\in J}\Vert U(s,x) - U(\bar s,x)\Vert_{L^2(\Omega)} & \le C_1\ q_1(|s-\bar s|),\notag\\
\sup_{t\in I}\Vert U(t,y) - U(t,\bar y)\Vert_{L^2(\Omega)} & \le C_2\ q_2(|y-\bar y|),
\end{align}
for some constants $C_1, C_2$. Then, arguing in a similar way as in Remark \ref{s2.2-r10}, we see that  $U$ satisfies \eqref{qoneqtwo}.
\end{rem}


For $\varepsilon\in(0,1)$, $j=(j_1,\ldots,j_{d_1}, j_{d_1+1},\ldots,j_{d_1+d_2})\in\mathbb{Z}^{d_1+d_2}$, define
\begin{align}
\label{2erres}
R^{\varepsilon,1}_j&=\prod_{i=1}^{d_1}\left[\frac{q_1^{-1}(\varepsilon)}{\sqrt{d_1}}j_i, \frac{q_1^{-1}(\varepsilon)}{\sqrt{d_1}}(j_i+1)\right],\ 
R^{\varepsilon,2}_j=\prod_{i=d_1+1}^{d_1+d_2}\left[\frac{q_2^{-1}(\varepsilon)}{\sqrt{d_2}}j_i, \frac{q_2^{-1}(\varepsilon)}{\sqrt{d_2}}(j_i+1)\right],\notag\\
 \tilde R_j^\varepsilon &= R^{\varepsilon,1}_j\times R^{\varepsilon,2}_j. 
\end{align}
For $t\in R^{\varepsilon,1}_j$ let $t_j^\varepsilon = \left(\frac{q_1^{-1}(\varepsilon)}{\sqrt{d_1}}j_i\right)_{i=1,\ldots,d_1}$, and for $x\in R^{\varepsilon,2}_j$, let $x_j^\varepsilon=\left(\frac{q_2^{-1}(\varepsilon)}{\sqrt{d_2}}j_i\right)_{i=d_1+1,\ldots,d_1+d_2}$.
\smallskip


Given two gauge functions $q_1$, $q_2$, and denoting $q=(q_1,q_2)$, we define
\beq
\label{gstbis}
\bar g_q(\tau) = \frac{\tau^D}{\left(q_1^{-1}(\tau)\right)^{d_1}  \left(q_2^{-1}(\tau)\right)^{d_2}},\quad \tau\in \R_+.
\eeq

\smallskip

\begin{teo}
\label{t2.2-tconclusive-st}
Let  $I$ and $J$ be compact subsets of $\R^{d_1}$ and $\R^{d_2}$, respectively, of positive Lebesgue measure. Assume that the stochastic process $U$ defined in \eqref{defU} is Gaussian, continuous, with i.i.d. components and such that $\sigma_{I,J}^2:=\inf_{(t,x)\in I^{(\eta)}\times J^{(\eta)}} {\rm Var}\ U(t,x)>0$ (for $\eta>0$ small enough). 
 Let $q_1$ and $q_2$ possess the same properties as the function $q$ in Lemma \ref{l2.2} and the function $\bar g_q$ given in \eqref{gstbis} be strictly increasing on a small interval $(0,\rho_0)$.
 Assume also that the process $U$ satisfies \eqref{qoneqtwo} on $I^{(\eta)}$ and $J^{(\eta)}$, respectively.
 
Then there exists a constant $C(I,J,\sigma_{I,J},d_1,d_2,D)$ such that for any Borel set $A\subset \R^D$, 
 \beq
\label{hgstbis}
P\left(U(I\times J)\cap A\ne \emptyset\right)\le C(I\times J,\sigma_{I,J},D,d_1,d_2)\mathcal{H}_{\bar g_q}(A).
\eeq
\end{teo}
\proof
Let $z\in\R^D$ and $\varepsilon>0$ (small enough). Since $U$ satisfies the conditions of Lemma \ref{l2.2-ts}, and  $\tilde R_j^\varepsilon \subset \tilde S_\varepsilon(t,x)$, 
for all $p\ge 1$, there exists a constant $C(p,d_1,d_2)$ such that for all $\varepsilon>0$ small enough,
\beq
\label{double-st-epsilon}
E\Big(\sup_{(t,x)\in \tilde R_j^\varepsilon}\left|U(t,x)-U(t_j^\varepsilon,x_j^\varepsilon)\right\vert^p\Big) \le C(p,d_1,d_2)\varepsilon^p. 
\eeq
Applying Theorem \ref{t2.3} with $M$ there replaced by $U$, we deduce
\beqn
P\left(U(\tilde R_j^\varepsilon)\cap B_\varepsilon(z)\neq \emptyset\right)
\le C(\sigma_{I,J},d_1,d_2,D) \varepsilon ^D,
\eeqn
for  some constant $C(\sigma_{I,J},d_1,d_2,D)$.
Similarly as in the proof of Lemma \ref{l2.3}, by an argument based on total probabilities (see \eqref{sbg-bis}), we deduce
\beq
\label{sbst-st}
P\left(U(I\times J)\cap B_\varepsilon(z)\ne \emptyset\right) \le C(I,J,\sigma_{I,J},d_1,d_2,D)\ \frac{\varepsilon^D}{\left(q_1^{-1}(\varepsilon)\right)^{d_1}\left(q_2^{-1}(\varepsilon)\right)^{d_2}}.
\eeq
We finish the proof by applying the covering argument laid out before the proof of Theorem \ref{2.2-tconclusive}.
\qed
\smallskip

Similarly as we did in Corollary \ref{s2-c10} and with the same arguments, from \eqref{sbst-st} we derive the following result on polarity of singletons.  
\begin{cor}
\label{s2-c10-bis}
The hypotheses are those of Theorem \ref{t2.2-tconclusive-st}. In addition assume that
\beq
\label{2.c-polar-bis}
\lim_{\tau\downarrow 0} \bar g_q(\tau) = 0.
\eeq
Then $P(U(I\times J)\cap \{z\}\ne\emptyset) =0$, that is, for the random field $U$ restricted to $I\times J$, any set $\{z\}\subset \R^D$ is polar.
\end{cor}
\smallskip

In Lemma \ref{s2.2-l-grow-g-st} (sections 1. and 2.) examples of gauge functions satisfying the conditions of Therorem \ref{2.2-tconclusive} and \eqref{2.c-polar-bis} are given.

\begin{rem}
\label{s2.1-r30}
The above discussion can be easily extended to the case where instead of $q_1$, $q_2$, we consider gauge functions $q_1, \ldots, q_{d_1+d_2}$ (repetitions are allowed). In this frame, setting $q=(q_1,\ldots,q_{d_1+d_2})$, and defining
\beq
\label{gstbis-multiple}
\bar g_q(\tau) = \frac{\tau^D}{\prod_{j=1}^{d_1+d_2}q_j^{-1}(\tau)},
\eeq
with suitable adaptation of conditions, we obtain \eqref{hgstbis} with $\bar g_q$ given in \eqref{gstbis-multiple}.This is \cite[Theorem 7.6, upper bound of (167), p. 188]{xia}.
\end{rem}
\subsection{Lower bounds for hitting probabilities}
\label{s2.1}
The aim of this section is to establish lower bounds on hitting probabilities of Gaussian processes in terms of 
$\textgoth{g}$-capacities.
\medskip

\noindent{\bf The single $q$-anisotropic case}
\medskip

Let  $M $ be a $D$-dimensional  Gaussian stochastic process, as given in \eqref{defM}, with i.i.d. components; let $K\subset \R^d$ be a compact set of positive Lebesgue measure.  We will use the notation $\sigma_x^2:={\text{Var}} (M_1(x))$, $\sigma_{x,\bar x}^2:= {\text{Cov}}(M_1(x),M_1(\bar x))$, 
$\rho_{x,\bar x} = {\text{Corr}}(M_1(x),M_1(\bar x))$, $f(x) = E(M(x))$, $\tilde M(x) = M(x) - f(x)$.
\smallskip

We introduce the following set of conditions.
\smallskip

\noindent{\em Hypotheses $(H_M)$ }  
\begin{enumerate}
\item There exist positive constants $c_1$, $c_2$ such that for all $x\in K$,
\beq
\label{2.1.1}
c_1\le \sigma_x^2\le c_2.
\eeq
\item $\rho_{x,\bar x} <1$ for all $x,\bar x\in K$.
\item There exist $\eta>0$ and $c_3>0$ such that for all $x,\bar x\in K$,
\beq
\label{2.1.2}
\left\vert \sigma_x^2 - \sigma_{\bar x}^2\right\vert \le c_3\  \Vert M_1(x)-M_1(\bar x)\Vert_{L^2(\Omega)}^{1+\eta}.
\eeq
\item There exists a gauge function $q$  such that for all $x,\bar x\in K$, 
\beq
\label{2.1.3}
\Vert \tilde M_1(x)-\tilde M_1(\bar x)\Vert_{L^2(\Omega)} \asymp q(|x-\bar x|), \quad |f(x) - f(\bar x)| \le C\ q(|x-\bar x|).
\eeq
\end{enumerate}
\medskip

\begin{rem}
\label{s2.1-r1}
Let ${\text{Var}}\, (M_1(\bar x)|M_1(x))$ denote the conditional covariance of $M_1(\bar x)$ given $M_1(x)$. The conditions 1 to 3 in  $(H_M)$ imply
\beq
\label{2.1.3-bis}
{\text{Var}}\, (M_1(\bar x)|M_1(x)) \asymp \Vert M_1(x)-M_1(\bar x)\Vert_{L^2(\Omega)}^2.
\eeq
Indeed, if $M$ is centred, this is \cite[Lemma 3.2, (1)]{dal:san} (with $\tau_{x,\bar x}^2:={\text{Var}}\ (M_1(\bar x)|M_1(x))$ there). Going through the proof we see that the property of being centred is not used. Along with \eqref{2.1.3}, we deduce
\beq
\label{2.1.3-tris}
{\text{Var}}\ (M_1(\bar x)|M_1(x)) = {\text{Var}}\ (\tilde M_1(\bar x) |\tilde M_1(x)) \asymp q^2(|x-\bar x|).
\eeq
\end{rem}
\medskip

Associated with the gauge function $q$ we define
\beq
\label{1.3-bis}
v_q(\tau) = \int_{q^{-1}(\tau)}^{\text{diam}(K)} [q(\rho)]^{-D} \rho^{d-1}\ d\rho,\ \tau\in\R_+.
\eeq
This section is devoted to prove the following statement.

\begin{teo}
\label{s2.1-t1}
Let $K\subset\R^d$ be a compact set of positive Lebesgue measure. Fix $N>0$ and let $A\subset B_N(0)\subset \R^D$ be a Borel set. 
Assume that conditions $(H_M)$ hold. Furthermore, suppose that 
\beq
\label{2.1.9}
\sup_{\tau\in[0, {\text{diam}}(K)]} v_q(\tau) g_q(\tau)\in (0,\infty),
\eeq
where $g_q$ is defined in \eqref{gg}. 
Then there exists a constant  $C:= C(f,K,N,d,D)>0$ such that
\beq
\label{2.1.10}
P\{M(K)\cap
 A \neq \emptyset\}\ge C\ {\text{Cap}}_{(g_q)^{-1}}(A).
\eeq

\end{teo}
\proof
We adapt the method used for example in \cite[Theorem 2.1]{b-l-x} inspired in 
\cite[pp. 204-206]{kahane}.

For any $x\in K$ and a probability measure $\mu$ on $A$, define
\begin{align}
\label{2.1.10-bis}
\bar \nu_n(x,\omega)&= \int_{A} (2\pi n)^{D/2} \exp\left(-\frac{n|M(x)- y|^2}{2}\right)\ \mu(dy)\notag\\
&=  \int_{A} \mu(dy) \int_{\R^D} d\xi  \exp\left(-\frac{|\xi|^2}{2n} + i\langle \xi, M(x)-y\rangle\right).
\end{align}
Consider the sequence of random measures on $K$, $(\nu_n, n\ge 1)$, with corresponding densities $(\bar \nu_n(x,\omega), n\ge 1)$. Set $\nu_n(K)(\omega) = \int_K \bar \nu_n(x,\omega)\ dx$. We aim to prove:
\begin{description}
\item {(i)} There exists $C_1>0$ such that for any $n\ge 1$,
$E\left( \nu_n(K)\right) \ge C_1$.
\item {(ii)} There exists $C_2>0$ such that for any $n\ge 1$,
$E\left[ \left(\nu_n(K)\right)^2\right] \le C_2\ \mathcal{E}_{(g_q)^{-1}}$.
\end{description}
By Paley-Zygmund inequality, this will imply
\beqn
P\{\nu_n(K) >0\} \ge \frac{\left[E\left( \nu_n(K)\right)\right]^2}{E\left[ \left(\nu_n(K)\right)^2\right]}\ge \frac{C_1}{C_2 \mathcal{E}_{(g_q)^{-1}}}.
\eeqn
Using an argument based on weak convergence of finite measures, we deduce \eqref{2.1.10}. 
\medskip

\noindent{\em Proof of (i).}  By Fubini's theorem,
\begin{align*}
E\left( \nu_n(K)\right)&= \int_K dx \int_{A} \mu(dy) \int_{\R^D} d\xi  \exp\left(-\frac{|\xi|^2}{2n} -i\langle \xi,y\rangle\right)
E\left(\exp(i\langle \xi,M(x)\rangle)\right)\\
&= \int_K dx \int_ A \mu(dy) \ \left(\frac{2\pi}{1/n+\sigma_x^2}\right)^{D/2}
\exp\left(-\frac{|y - f(x)|^2}{2[1/n+\sigma_x^2]}\right).
\end{align*}
The last equality is obtained computing first the characteristic function \break $E\left(\exp(i\langle \xi,M(x)\rangle)\right)$ 
and then, the Fourier inversion formula.

Let $N_0 = N+\sup_{x\in K}|f(x)|$. Applying \eqref{2.1.1}, and since on the set $A$, $|y-f(x)|\le N_0$, the above computations yield
\begin{align*}
E\left( \nu_n(K)\right)&\ge  \int_K dx \int_ A \mu(dy)\left(\frac{2\pi}{1+\sigma_x^2}\right)^{D/2}
\exp\left(-\frac{N_0^2}{2\sigma_x^2}\right)\\
&\ge |K|\left(\frac{2\pi}{1+c_2}\right)^{D/2} \exp\left(-\frac{N_0^2}{2c_1}\right):=C_1.
\end{align*}
This ends the proof of (i). Notice that $C_1:=C_1(f,K,N,D)$.
\medskip

\noindent{\em Proof of (ii).}  For any $x, \bar x\in K$, $y, \bar y\in A$, set
\begin{align*}
&I(x,\bar x, y, \bar y)\\
&\quad = \int_{\R^D \times \R^D} e^{-i\langle (\xi,\bar \xi), (y,\bar y)\rangle}
\exp\left(-\frac{|(\xi,\bar\xi)|^2}{2n}\right)
\exp\left(i\langle(\xi,\bar\xi), (M(x),M(\bar x))\rangle\right)d\xi\ d\bar\xi.
\end{align*}
Using \eqref{2.1.10-bis}, the definition of $\nu_n(K)$ and Fubini's theorem, we see that
\beq
\label{2.1.13}
E\left[ \left(\nu_n(K)\right)^2\right] = \int_{K\times K} dx\ d\bar x\int_{A\times A} \mu(dy)\ \mu(d\bar y) 
E\left(I(x,\bar x, y, \bar y)\right).
\eeq
With elementary computations based on the properties of the exponential function, we deduce $I(x,\bar x, y, \bar y)=\prod_{j=1}^D\  I_j(x,\bar x, y, \bar y)$, 
with
\begin{align*}
 &I_j(x,\bar x, y, \bar y)\\
 &\ =\int_{\R^2} d\xi_j\ d\bar\xi_j\
e^{-i\langle(\xi_j,\bar\xi_j), (y_j,\bar y_j)\rangle} 
\exp\left(-\frac{|(\xi_j,\bar\xi_j)|^2}{2n}\right)
\exp\left(i\langle(\xi_j,\bar\xi_j),(M_j(x),M_j(\bar x))\rangle\right)
\end{align*}
Since the factors in the product above are i.i.d random variables, from \eqref{2.1.13} we obtain
\beq
\label{2.1.14}
E\left[ \left(\nu_n(K)\right)^2\right] = \int_{K\times K} dx\ d\bar x\int_{A\times A} \mu(dy)\ \mu(d\bar y) 
\prod_{j=1}^D \left[E\left(I_j(x,\bar x, y, \bar y)\right)\right].
\eeq
Let $\Gamma_{x,\bar x}$ denote the covariance matrix of the $2$-dimensional Gaussian random vector $(M_j(x),M_j(\bar x))$ (which is the same as for $(\tilde M_j(x), \tilde M _j(\bar x))$),
and set $\Gamma_{x,\bar x}^n = \frac{1}{n} {\text{Id}}_2 + \Gamma_{x,\bar x}$. 
Computing 
$E\left(\exp\left(i\langle(\xi_j,\bar\xi_j),(M_j(x),M_j(\bar x))\rangle\right)\right)$ and then applying the Fourier inversion formula, we obtain
\begin{align}
\label{2.1.15}
&E\left(I_j(x,\bar x, y, \bar y)\right)\\
&\  = \int_{\R^2} d\xi_j\ d\bar\xi_j\
e^{-i\langle(\xi_j,\bar\xi_j), (y_j-f_j(x),\bar y_j-f_j(\bar x))\rangle} 
\exp\left(-\frac{1}{2}(\xi_j,\bar\xi_j) \Gamma_{x,\bar x}^n(\xi_j,\bar\xi_j)^{\intercal}\right)\notag\\
&\ = \frac{2\pi}{\left(\det\Gamma_{x,\bar x}^n\right)^{1/2}}
\exp\left(-\frac{1}{2} (y_j-f_j(x), \bar y_j-f_j(\bar x))\left(\Gamma_{x,\bar x}^n\right)^{-1}(y_j-f_j(x), \bar y_j-f_j(\bar x))^{\intercal}\right).
\end{align} 
Explicit computations show
\begin{align*}
(y_j-f_j(x), \bar y_j-f_j(\bar x))\left(\Gamma_{x,\bar x}^n\right)^{-1}&(y_j-f_j(x), \bar y_j-f_j(\bar x))^{\intercal}\\
&\geq 
\frac{E\left[\left((y_j-f_j(x)) \tilde M_j(x) - (\bar y_j -f_j(\bar x))\tilde M_j(\bar x)\right)^2\right]}{\det \Gamma^n_{x,\bar x}}.
\end{align*}
Hence, applying Lemma \ref{s2.1-l1} we deduce
\beq
\label{2.1.16}
E\left(I_j(x,\bar x, y, \bar y)\right)  \le C \frac{1}{\left(\det \Gamma^n_{x,\bar x}\right)^{1/2}}
\exp\left(-\frac{c|(y_j-\bar y_j)-(f_j(x)-f_j(\bar x))|^2}{2 \det \Gamma^n_{x,\bar x}}\right).
\eeq
Using this estimate in \eqref{2.1.14}, we obtain
\begin{align}
\label{2.1.17-bbis}
E\left[ \left(\nu_n(K)\right)^2\right]
&\le C \int_{K\times K} dx\ d\bar x\int_{A\times A} \mu(dy)\ \mu(d\bar y) 
\frac{1}{\left(\det \Gamma^n_{x,\bar x}\right)^{D/2}}\notag\\
&\qquad \qquad \quad\times\exp\left(-\frac{c|(y-\bar y)-(f(x)-f(\bar x))|^2}{2 \det \Gamma^n_{x,\bar x}}\right).
\end{align}

 Since $\Gamma_{x,\bar x}$ is nonnegative definite,
  \beq
  \label{detlow}
 \det \Gamma_{x,\bar x}^n\ge \det \Gamma_{x,\bar x} = \sigma_x^2 {\text{Var}}\ (M(\bar x)|M(x)) \ge C q^2(|x-\bar x|),
 \eeq
 where the last inequality follows from \eqref{2.1.1} and \eqref{2.1.3-tris}. This estimate along with \eqref{2.1.3} implies
  \beq
 \label{exp-bound}
 \sup_{x,\bar x \in K}\frac{|f(x)-f(\bar x)|^2}{\det \Gamma^n_{x,\bar x}}\le C<\infty.
 \eeq
 Apply the inequality $|(y-\bar y)-(f(x)-f(\bar x))|^2\ge \frac{1}{2}|y-\bar y|^2 - |f(x)-f(\bar x)|^2$ and \eqref{exp-bound} on the right-hand side of \eqref{2.1.17-bbis} to deduce,
\beq
\label{2.1.17}
E\left[ \left(\nu_n(K)\right)^2\right]
 \le C \int_{K\times K} dx\ d\bar x\int_{A\times A} \mu(dy)\ \mu(d\bar y) 
\frac{1}{\left(\det \Gamma^n_{x,\bar x}\right)^{D/2}}
\exp\left(-\frac{c|y-\bar y|^2}{2 \det \Gamma^n_{x,\bar x}}\right).
\eeq

If $\det \Gamma^n_{x,\bar x}\ge |y-\bar y|^2$, the integrand is bounded from above by the factor $\left(\det \Gamma^n_{x,\bar x}\right)^{-D/2}$. If on the contrary, $\det \Gamma^n_{x,\bar x}< |y-\bar y|^2$, the integrand is bounded (up to a multiplicative constant) by $|y-\bar y|^{-D}$, because  the function
$z\mapsto z^{D/2} e^{-cz}$ is bounded over $\R_+$. In this way,
\begin{align}
\label{2.1.18}
E\left[ \left(\nu_n(K)\right)^2\right] &\le C\int_{K\times K} dx\ d\bar x\int_{A\times A} \mu(dy)\ \mu(d\bar y) 
\frac{1}{\max\left( \left( \det \Gamma_{x,\bar x}^n\right)^{D/2}, |y-\bar y|^D\right)}\notag\\
&\le C\int_{K\times K} dx\ d\bar x\int_{A\times A} \mu(dy)\ \mu(d\bar y)\  
\frac{1}{\max\left(q^D(|x-\bar x|), |y-\bar y|^D\right)},
\end{align}
where in the second inequality we have applied \eqref{detlow}.
 
 Our next goal is to prove 
 \beq
 \label{2.1.19}
 \int_{K\times K}\ \frac {dx\ d\bar x}{\max\left(q^D(|x-\bar x|), |y-\bar y|^D\right)} \le C(K,d) [g_q(|y-\bar y|)]^{-1}.
 \eeq
 Indeed, 
 \begin{align*}
& \int_{(K\times K)\cap \{q(|x-\bar x|)\le |y-\bar y|\}} \frac{ dx\ d\bar x}{\max\left(q^D(|x-\bar x|), |y-\bar y|^D\right)}\\  
& \qquad\qquad\le C(K,d)|y-\bar y|^{-D} \int_0^{q^{-1}(|y-\bar y|)} \rho^{d-1} d\rho = C(K,d)\  [g_q(|y-\bar y|)]^{-1},
  \end{align*}
  and
   \begin{align*}
 &\int_{(K\times K)\cap \{q(|x-\bar x|)>|y-\bar y|\}} \frac{dx\ d\bar x}{\max\left(q^D(|x-\bar x|), |y-\bar y|^D\right)}\\
 &\ \le  C(K,d) \int_{q^{-1}(|y-\bar y|)}^{\text{diam}\ (K)} [q(\rho)]^{-D} \rho^{d-1} d\rho  = C(K,d) v_q(|y-\bar y|)\le \tilde C(K,d)  [g_q(|y-\bar y|)]^{-1},
 \end{align*}
 where the last equality holds because of hypothesis \eqref{2.1.9}.
 
 Hence,
 \beq
 \label{2.1.20}
 E\left[ \left(\nu_n(K)\right)^2\right]  \le C(K,d) \mathcal{E}_{(g_q)^{-1}}(\mu),
 \eeq
 and the right-hand side does not depend of $n$.
 
 The proof of the theorem is complete.
 \qed

Consider  the gauge functions $q$ introduced in Examples \ref{2.2-e0}. In Lemma \ref{s2.1-l2} we compute the corresponding functions
 $v_q$ (defined in \eqref{1.3-bis}) and find $g_q$ satisfying  condition \eqref{2.1.9}.
\smallskip

We end this section with a technical result used in the proof of Theorem \ref{s2.1-t1}.
\begin{lem}
\label{s2.1-l1}
Assume $(H_M)$. Then for any $a, b \in \R$ and $x, \bar x\in K$, there exists a constant $c>0$ such that
\beq
\label{2.1.30}
E[(a \tilde M_1(x) - b \tilde M_1(\bar x))^2] \ge c(a-b)^2.
\eeq
\end{lem}
\proof
Property \eqref{2.1.30} is equivalent to say that the matrix
\beqn
N_{x,\bar x}
\begin{pmatrix}
\sigma_x^2 - c & -(\sigma^2_{x,\bar x} - c)\\
 -(\sigma^2_{x,\bar x} - c) & \sigma_{\bar x}^2-c
 \end{pmatrix}
 \eeqn
is nonnegative definite. Computing $\det N_{x,\bar x}$, we see that this holds if and only if \break
$
\left(\det \Gamma_{x,\bar x}\right)\left(\Vert\tilde M_1(x)  - \tilde M_1(\bar x)\Vert_{L^2(\Omega)}^2\right)^{-1}\ge c$.
Applying Remark \ref{s2.1-r1}, and using \eqref{2.1.1}, we obtain
\beqn
\frac{\det \Gamma_{x,\bar x}}{\Vert\tilde M_1(x)  - \tilde M_1(\bar x)\Vert_{L^2(\Omega)}^2}
= \frac{\sigma_x^2 {\text{Var}} (\tilde M_1(\bar x)\big |\tilde M_1(x))}{\Vert\tilde M_1(x)  - \tilde M_1(\bar x)\Vert_{L^2(\Omega)}^2}
\asymp 1.
\eeqn
\qed

\bigskip


\noindent{\bf The multiple $q$-anisotropic case}
\medskip

Let $U=\{U(t,x)=(U_1(t,x), \ldots, U_D(t,x)),\ (t,x)\in\R^{d_1}\times \R^{d_2}\}$ be a $D$-dimensional Gaussian stochastic process 
with i.i.d. components, and $I\subset\R^{d_1}$, $J\subset \R^{d_2}$ be compact sets of positive Lebesgue measure.  We will use the notation $\sigma_{t,x}^2:={\text{Var}} (U_1(t,x))$, $\sigma_{(t,x),(s,y)}^2:= {\text{Cov}}(U_1(t,x),U_1(s,y))$, $\rho_{(t,x),(s,y)} = {\text{Corr}}(U_1(t,x),U_1(s,y))$, $f(t,x) = E(U(t,x))$ and $\tilde U(t,x) = U(t,x) - f(t,x)$. 
\smallskip

By analogy with assumptions  $(H_M)$ in the discussion on single $q$-anisotropic processes, we introduce
the following set of conditions.
\smallskip

\noindent{\em Hypotheses $(H_U)$}
\begin{enumerate}
\item There exist positive constants $c_1$, $c_2$ such that for all $(t,x)\in I\times J$,
\beq
\label{2.1.1-bis}
c_1\le \sigma_{t,x}^2\le c_2.
\eeq
\item $\rho_{(t,x),(s,y)} <1$ for all $(t,x), (s,y)\in I\times J$.
\item There exist $\eta>0$ and $c_3>0$ such that for all $(t,x), (s,y)\in I\times J$,
\beq
\label{2.1.2-bis}
\left\vert \sigma_{t,x}^2 - \sigma_{s,y}^2\right\vert \le c_3\  \Vert U_1(t,x)-U_1(s,y)\Vert_{L^2(\Omega)}^{1+\eta}.
\eeq
\item There exist  gauge functions $q_1$, $q_2$ such that for all $(t, x), (s,y)\in I\times J$, 
\begin{align}
\Vert \tilde U_1(t,x)-\tilde U_1(s,y)\Vert_{L^2(\Omega)} &\asymp q_1(|t-s|) + q_2(|x-y|). \label{2.1.3-bisbis}\\
|f(t,x) - f(s,y)| &\le C \left(q_1(|t-s|) + q_2(|x-y|)\right). \label{upexpect}
\end{align}\end{enumerate}
\medskip

\begin{rem}
\label{s2.1-r2}
Using similar arguments as in Remark \ref{s2.1-r1}, now applied to the process $U$, we deduce
\beq
\label{2.1.50}
\text{Var}\,  \left(U_1(t,x)| U_1(s,y)\right)\asymp \Vert U_1(t,x)-U_1(s,y)\Vert_{L^2(\Omega)}^2 \asymp \left(q_1(|t-s|) + q_2(|x-y|)\right)^2,
\eeq
for any $(t, x), (s,y)\in I\times J$.
\end{rem}
Set $d_I = \text{diam}\,(I)$, $d_J= \text{diam}\,(J)$, $c_{I,J} = \max(q_1(d_I), q_2(d_J))$. Assuming that $q_1$ and $q_2$ are differentiable,  for $\tau\in[0, c_{I,J}]$, define
\beq
\label{2.1.51}
\bar v_q(\tau)= \int_\tau^{c_{I,J}} \rho^{-D+1} \left[q_1^{-1}(\rho)\right]^{d_1-1}\left[q_2^{-1}(\rho)\right]^{d_2-1}
\left[\dot q_1(q_1^{-1}(\rho))\right]^{-1} \left[\dot q_2(q_2^{-1}(\rho))\right]^{-1}\ d\rho.
\eeq

To highlight the analogy between $\bar v_q$ and the function $v_q$ defined in \eqref{1.3-bis}, we observe that if $q$ in \eqref{1.3-bis} is differentiable, with the change of variable $\rho\mapsto q(\rho)$ we have
\beqn
v_q(\tau) = \int_{\tau}^{q(\text{diam}(K))} \rho^{-D}  \left[q^{-1}(\rho)\right]^{d-1}\left[\dot q(q^{-1}(\rho))\right]^{-1} d\rho.
\eeqn
\smallskip

Our purpose is to prove the following result.

\begin{teo}
\label{s2.1-t2}
Let  $I\subset\R^{d_1}$ and $J\subset \R^{d_2}$ be compact sets of positive Lebesgue measure. 
Fix $N>0$ and let $A\subset B_N(0)\subset \R^D$ be a Borel set. 
Assume that conditions $(H_U)$ hold. Furthermore,  
suppose that on $(0,c_{I,J})$, the gauge functions $q_i$, $i=1,2$, are differentiable with decreasing derivatives  $\dot q_i$, and
\beq
\label{2.1.52}
\sup_{\tau\in[0, c_{I,J}]} \bar v_q(\tau/2) \bar g_q(\tau)\in (0,\infty),
\eeq
where $\bar g_q$ is the function defined in \eqref{gstbis}.

Then there exists a constant  $C:= C(f,I,J,N,d_1,d_2,D)>0$ such that
\beq
\label{2.1.10-bis-bis}
P\{U(I\times J)\cap
 A \neq \emptyset\}\ge C\ {\text{Cap}}_{(\bar g_q)^{-1}}(A).
\eeq
\end{teo}

\proof
The approach to the proof is the same as that of Theorem \ref{s2.1-t1}. To avoid repetitions, we only provide details on the relevant differences. 

For any $(t,x)\in I\times J$ and a probability measure $\mu$ on $A$, define 
\begin{align}
\label{2.1.10-bisbis}
\bar \nu_n((t,x),\omega)&= \int_{A} (2\pi n)^{D/2} \exp\left(-\frac{n|U(t,x)- y|^2}{2}\right)\ \mu(dy)\notag\\
&=  \int_{A} \mu(dy) \int_{\R^D} d\xi  \exp\left(-\frac{|\xi|^2}{2n} + i\langle \xi, U(t,x)-y\rangle\right),\quad n\ge 1,
\end{align}
and let $\nu_n(I\times J)(\omega) = \int_{I\times J}\bar \nu_n((t,x),\omega)\ dt\,dx$.

Applying \eqref{2.1.1-bis}, similarly as for the proof of (i) in Theorem \ref{s2.1-t1}, we obtain
\beq
\label{2.1.11-bis}
E\left(\nu_n(I\times J)\right) \ge \bar C_1,
\eeq
with $\bar C_1=C_1(f,I,J,N,D)$.

With similar computations as those used to derive \eqref{2.1.18}, we have
\begin{align}
\label{2.1.18-bis}
&E\left[ \left(\nu_n(I\times J)\right)^2\right]\notag\\ 
& \le C\int_{(I\times J)^2} dt\,dx\,d\bar t\,d\bar x\int_{A\times A} \mu(dy)\ \mu(d\bar y)\  
\frac{1}{\max\left([q_1(|t-\bar t|)+q_2(|x-\bar x|)]^D, |y-\bar y|^D\right)}.
\end{align}
For $h\ge 0$, set
\beq
\label{2.1.60}
I:= \int_{(I\times J)^2} dt\,dx\,d\bar t\,d\bar x \left[\max\left([q_1(|t-\bar t|)+q_2(|x-\bar x|)]^D,h^D\right)\right]^{-1}.
\eeq
Apply the change of variables $(t,\bar t)\mapsto (t, t-\bar t)$, $(x,\bar x)\mapsto (x, x-\bar x)$, to deduce
\beq
\label{2.1.61}
I\le |I\times J| \int _{B_{d_I}(0)} dr  \int _{B_{d_J}(0)} dz \left[\max\left([q_1(|r|)+q_2(|z|)]^D,h^D\right)\right]^{-1},
 \eeq
 where $|I\times J|$ denotes the Lebesgue measure of $I\times J$.
  
 Let $I_1$ denote the integral in \eqref{2.1.61} over the set of points $(r,z)$ satisfying $q_1(|r|)+q_2(|z|)\le h$. 
 Changing to polar coordinates, we see that
 \begin{align}
 \label{2.1.62}
 I_1 &= h^{-D} \int _{B_{d_I}(0)} dr  \int _{B_{d_J}(0)} dz\ 1_{\{q_1(|r|)+q_2(|z|)\le h\}}\notag\\
 &\le h^{-D} \left(\int _{B_{d_I}(0)} dr\,  1_{\{q_1(|r|)\le h\}}\right)
 \left(\int _{B_{d_J}(0)} dz\ 1_{\{q_2(|z|)\le h\}}\right)\notag\\
 &\le  C(d_1,d_2) h^{-D}\left(\int_0^{q_1^{-1}(h)} \rho^{d_1-1} d\rho\right)
 \left(\int_0^{q_2^{-1}(h)} \rho^{d_2-1} d\rho\right)
 = C(d_1,d_2) \left[\bar g_q(h)\right]^{-1}.
 \end{align}
 
 Next, we denote by $I_2$ the integral in \eqref{2.1.61} over the set of points $(r,z)$ such that $q_1(|r|)+q_2(|z|)> h$. Applying two changes of variables: first polar coordinates, $r\mapsto(\rho_1,\theta_1)$, $z\mapsto(\rho_2,\theta_2)$, and then $\rho_i\mapsto q_i(\rho_i)$, $i=1,2$, we obtain
 \begin{align*}
 I_2&= C(d_1,d_2) \int_0^{d_I} d\rho_1 \int_0^{d_J}d\rho_2\ 1_{\{q_1(\rho_1)+q_2(\rho_2)> h\}}
\left[q_1(\rho_1) + q_2(\rho_2)\right]^{-D}  \rho_1^{d_1-1}   \rho_2^{d_2-1}\\
 &= C(d_1,d_2) \int_0^{q_1(d_I)} d\tau_1\int_0^{q_2(d_J)} d\tau_2\ 1_{\{\tau_1+\tau_2>h\}}(\tau+\tau_2)^{-D}
 \\
 &\qquad \times 
 \left(q_1^{-1}(\tau_1)\right)^{d_1-1} \left(q_2^{-1}(\tau_2)\right)^{d_2-1}\left[\dot q_1(q_1^{-1}(\tau_1))\right]^{-1} \left[\dot q_2(q_2^{-1}(\tau_2))\right]^{-1}\\
 &\le C(d_1,d_2, D)  \int_0^{q_1(d_I)} d\tau_1\int_0^{q_2(d_J)} d\tau_2\ 1_{\{|(\tau_1,\tau_2)|>h/2\}}
[ |(\tau_1,\tau_2)|]^{-D}\\
&\qquad \times \left(q_1^{-1}(\tau_1)\right)^{d_1-1} \left(q_2^{-1}(\tau_2)\right)^{d_2-1}
 \left[\dot q_1(q_1^{-1}(\tau_1))\right]^{-1} \left[\dot q_2(q_2^{-1}(\tau_2))\right]^{-1},\\
 \end{align*}
where in the last inequality we have used  $|(\tau_1,\tau_2)|\le \tau_1+\tau_2\le 2 |(\tau_1,\tau_2)|$ ($|\cdot|$ is the Euclidean norm). Changing $(\tau_1,\tau_2)$ into polar coordinates, because for $i=1,2,$ $q_i$ are increasing and $\dot q_i$ decreasing, we deduce
 \begin{align}
 \label{2.1.63}
 I_2 & \le  C(d_1,d_2,D)\int_{h/2}^{c_{I,J}} \rho^{-D+1} \left[q_1^{-1}(\rho)\right]^{d_1-1}\left[q_2^{-1}(\rho)\right]^{d_2-1}\notag\\
&\qquad \qquad\qquad\times\left[\dot q_1(q_1^{-1}(\rho))\right]^{-1} \left[\dot q_2(q_2^{-1}(\rho))\right]^{-1}\ d\rho\notag\\
&  = C(d_1,d_2,D) \bar v_q(h/2) \le C(I,J) \left[\bar g_q(h)\right]^{-1},
\end{align}
where the last inequality follows from the assumption \eqref{2.1.52}.
 
Thus, from  \eqref{2.1.18-bis} by applying \eqref{2.1.62} and \eqref{2.1.63} with $h:=|y-\bar y|$, we obtain
\beq
\label{2.1.20-bis}
E\left[\left(\nu_n(I\times K)\right)^2\right] \le C(I,J,d_1,d_2,D)\, \mathcal{E}_{(\bar g_q)^{-1}} (\mu).
\eeq
We conclude in a similar way as in the proof of Theorem \ref{s2.1-t1}.
\qed

In Lemma \ref{s2.1-l2}, we give two examples where Theorem \ref{s2.1-t2} can be applied. 
\medskip

By the definition of capacity, we have (see e.g. \cite[p.529]{kho})
\beq
\label{s2.1.65}
\text{Cap}_{\textgoth{g}}(\{z\}) >0\quad  \text{if and only if} \quad \textgoth{g}(0) <\infty.
\eeq
Take $A=\{z\}$, $z\in \R^D$, in Theorems \ref{s2.1-t1} and \ref{s2.1-t2}. If $\{z\}$ is polar for the  process $M$ restricted to the compact $K$ (respectively, for the process $U$ restricted to the compact $I\times J$), then necessarily, $\text{Cap}_{(g_q)^{-1}}(\{z\})=0$ (respectively,  $\text{Cap}_{(\bar g_q)^{-1}}(\{z\})=0$. According to \eqref{s2.1.65} this is equivalent to $g_q(0)=0$ (respectively,  $\bar g_q(0)=0$). Together with Corollaries \ref{s2-c10} and \ref{s2-c10-bis} we obtain the following result on polarity of points.

\begin{prop}
\label{s2.1-p-polar}
A singleton $\{z\}$ is polar for the process $M$ restricted to the compact $K$ (respectively, for the process $U$ restricted to the compact $I\times J$) if and only if $\lim_{\tau\downarrow 0} g_q(\tau)=0$ (respectively,  $\lim_{\tau\downarrow 0} \bar g_q(\tau)=0$).
\end{prop}

\subsection{Examples}
\label{s2.3}

Under the unifying umbrella provided by Theorems \ref{2.2-tconclusive}, \ref{t2.2-tconclusive-st}, \ref{s2.1-t1} and \ref{s2.1-t2}, we present in this section a selection of known results on hitting probabilities. We defer to Section \ref{s3} the new application to the multiple $q$-anisotropic process that has motivated this work.

\begin{exa}
\label{2.2-e2}
Fix compact sets $I\subset\R^{d_1}$, $J\subset\R^{d_2}$ of positive Lebesgue measure, and $\varepsilon\in(0,1)$. Assume that the process $U$ defined in \eqref{defU} is Gaussian with i.i.d. components. Suppose that there exist $\nu_1, \nu_2\in(0,1)$ and for any $(s,y), (t,x)\in(I\times J)^{(2\varepsilon)}$, 
\beq
\label{hgst-example-bis-assump}
\left\Vert U(t,x)-U(s,y)\right\Vert_{L^2(\Omega)} \asymp \left( |t-s|^{\nu_1} + |x-y|^{\nu_2}\right).
\eeq
By Kolmogorov's continuity theorem, $\{U(t,x)\}_{(t,x)\in I\times J}$ has continuous sample paths, a.s.

The condition \eqref{qoneqtwo} holds with $q_i(\tau)= \tau^{\nu_i}$, $i=1,2$ (see Remark \ref{s2.2-r20}). The function $\bar g_q$ defined in \eqref{gstbis} is
\beqn
\bar g_q(\tau)= r^{D-\left(\frac{d_1}{\nu_1}+\frac{d_2}{\nu_2}\right)},
\eeqn
and it is increasing if $D>(\tfrac{d_1}{\nu_1}+\tfrac{d_2}{\nu_2})$ (see Lemma \ref{s2.2-l-grow-g-st}). 

Assume that the process $U$ satisfies $\sigma_{I,J}^2:=\inf_{(t,x)\in(I\times J)^{(2\varepsilon)}} {\rm Var}\ (U(t,x))>0$; then from Theorem \ref{t2.2-tconclusive-st} we deduce the following: 

There exists a constant $C:=C(I\times J,\sigma_{I,J},D,d_1,d_2)$ such that 
for any Borel set $A\subset\R^D$, 
\beq
\label{hgst-example-bis}
P\left(U(I\times J)\cap A\ne \emptyset\right)\le C\,\mathcal{H}_{D-\left(\frac{d_1}{\nu_1}+\frac{d_2}{\nu_2}\right)}(A).
\eeq
If $D\le (\tfrac{d_1}{\nu_1}+\tfrac{d_2}{\nu_2})$, by definition of the Hausdorff measure, $\mathcal{H}_{D-\left(\frac{d_1}{\nu_1}+\frac{d_2}{\nu_2}\right)}(A)=\infty$. Thus \eqref{hgst-example-bis}  is still valid, but non informative.
\smallskip

By  Lemma \ref{s2.1-l2} (1.) we deduce the validity of \eqref{2.1.52}. Therefore, assuming hat $U$ satisfies   $(H_U)$, we see that the hypotheses of Theorem \ref{s2.1-t2} are satisfied. Thus, 
for any bounded Borel set $A\subset B_N(0)\subset \R^D$
there exists  $c:= c(I,J,d_1,d_2,D)$  such that
\beq
\label{2.1.40}
P\{U(I\times J)\cap A \neq \emptyset\}\ge c\, {\text{Cap}}_{D-\left(\tfrac{d_1}{\nu_1}+\tfrac{d_2}{\nu_2}\right)}(A)  
 .
\eeq

With \eqref{hgst-example-bis} and  \eqref{2.1.40}, we recover a version of \cite[Theorem 7.6]{xia} on hitting probabilities in the {\em classical} centred anisotropic case. 

Remark \ref{s2.1-r30} motivates an extension of \eqref{hgst-example-bis}. Indeed, let us replace the upper bound in 
\eqref{hgst-example-bis-assump} by 
\beq
\label{hgst-example-bis-assump-bis}
\left\Vert U(t,x)-U(s,y)\right\Vert_{L^2(\Omega)} \le C\left( \sum_{j=1}^{d_1}|t_j-s_j|^{\delta_j} +\sum_{k=1}^{d_2} |x_k-y_k|^{\nu_k}\right).
\eeq
Then, taking $q_j(\tau)= \tau^{\delta_j}$, $j=1,\ldots, d_1$, and $q_j(\tau) = \tau^{\nu_j}$, $j=d_1+1,\ldots, d_1+d_2$, we have
\beqn
\bar g_q(\tau) = \frac{\tau^D}{\prod_{j=1}^{d_1} \tau^{1/\delta_j} \prod_{k=1}^{d_2} \tau^{1/\nu_k}}.
\eeqn
Consequently, letting 
$Q:= \left(\sum_{j=1}^{d_1}\tfrac{1}{\delta_j}\right)+ \left(\sum_{k=1}^{d_2}+\tfrac{1}{\nu_k}\right)$,
we obtain
\beq
\label{hgst-example-tris}
P\left(U(I\times J)\cap A\ne \emptyset\right)\le C(I\times J,\sigma_{I,J},D,d_1,d_2)\mathcal{H}_{D-Q}(A).
\eeq

The following examples provide illustrations of the preceding results. 
\smallskip

{\em 1.}\ Consider the random field solution, $\{u(t,x), (t,x)\in[0,T]\times [0,L]\}$, to the system of linear stochastic heat equations
\beqn
\left(\frac{\partial}{\partial t}- \frac{\partial^2}{\partial x^2}\right) u_i(t,x) = \dot W_i(t,x), \quad  (t,x)\in(0,T]\times[0,L],\ i=1,\ldots,D,
\eeqn
with null initial value and vanishing Dirichlet boundary conditions. The processes \break
$(\dot W_i(t,x))_i$ are independent space-time white noises. Here, $d_1=d_2=1$, $\nu_1=1/4$ and $\nu_2= 1/2$ (see e.g. \cite{dal:kho}). Hence, \eqref{hgst-example-bis} and \eqref {2.1.40} hold with $\mathcal{H}_{D-6}(A)$ and
${\text{Cap}}_{D-6}(A)$, respectively. We recover  
\cite[Theorems  2.1 and  3.1]{dal:kho}. 
\smallskip

{\em 2.}\  For any $k\ge 1$, let $\{u(t,x), (t,x)\in(0,T]\times \R^k\}$ be the random field solution to the system of linear stochastic wave equations
\beqn
\left(\frac{\partial^2}{\partial t^2}-\frac{\partial^2}{\partial x^2}\right) u_i(t,x) = \dot W_i(t,x),\quad  (t,x)\in(0,T]\times \R^k,\  i=1,\ldots,D,
\eeqn
with null initial conditions. Assume that $(\dot W_i(t,x))_i$  are independent 
 noises white in time, with a stationary spatial covariance given by a Riesz kernel of order $\beta\in(0, k\wedge 2)$. In this example, $d=1+k$. According to the results in \cite{dal:san}, on any set $I\times J=[t_0,T]\times [-M,M]^k$  $(t_0, M>0)$, the hypotheses introduced above are satisfied with $\nu_1=\nu_2 = \frac{2-\beta}{2}$ in \eqref{hgst-example-bis-assump}. Therefore, 
 $\mathcal{H}_{g_q}(A) = \mathcal{H}_{D-2(1+k)/(2-\beta)}(A)$ and $\text{Cap}_{g_q}(A) = \text{Cap}_{D-2(1+k)/(2-\beta)}(A)$. We therefore recover  \cite[Theorems 4.4 and 4.5]{dal:san}.
 \end{exa}

\begin{exa}
\label{2.2-e-q}
Fix a compact set $K\subset \R^d$ of positive Lebesgue measure, and $\varepsilon\in(0,1)$. Let  the process $M$ in \eqref{defM}
 be Gaussian, centred, with i.i.d. components. Suppose there exists $\nu\in(0,1)$ such that,  for any $x,y\in K^{(2\varepsilon)}$,
 \beq
 \label{canon}
\left\Vert M(x) - M(y)\right\Vert_{L^2(\Omega)}\asymp |x-y|^\nu.
\eeq
By Kolmogorov's continuity theorem, the sample paths of $\{M(x)\}_{x\in K}$ are continuous a.s.
Appealing to Remark \ref{s2.2-r10}, we see that \eqref{eq2.2.17} holds with $q(\tau)= \tau^{\nu}$ and hence,
$g_q(\tau)= \tau^{D-\frac{d}{\nu}}$. Assume that $\nu>d/D$ (which, according to Lemma \ref{s2.2-l-grow-g-st} (1.), ensures that $g_q$ is increasing) and furthermore, $\sigma_K^2:=\inf_{x\in K^{(\eta)}} {\rm{Var}}\ (M(x))>0$ (for $\eta>0$ small enough). From Theorem \ref{2.2-tconclusive}, we deduce, for any Borel set $A\subset\R^D$, 
\beq
\label{hgst-example}
P\left(M(K)\cap A\ne \emptyset\right)\le C\, \mathcal{H}_{D-\frac{d}{\nu}}(A),
\eeq
with $C:= C(K,\sigma_K,d,D)$. If $\eta\le d/D$, the right-hand side in \eqref{hgst-example} is infinite (by definition). Thus  \eqref{hgst-example} still holds.

In addition to the above assumptions, suppose that $M$ satisfies $(H_M)$ and observe that by Lemma \ref{s2.1-l2} (1.), \eqref{2.1.9} holds.  Fix a bounded Borel set $A\subset B_N(0)\subset\R^D$. Then applying Theorem \ref{s2.1-t1} we obtain,
\beq
\label{hgst-example-lb}
P\left(M(K)\cap A\ne \emptyset\right)\ge c\, \text{Cap}_{D-\frac{d}{\nu}}(A),
\eeq
with $c:= c(K,N,d,D)$.

Consider the example of a system of linear stochastic Poisson equations on an open set $O\subset \R^d$  $(d=1,2,3)$:
 \beq
 \label{ex-poisson}
 -\Delta u_i(x) = \sum_{j=1}^d \sigma_{i,j}\dot W^j(x),\ x\in O, \ i=1,\ldots,D,\qquad u_{|\partial O}=0,
 \eeq
where $\dot W = (\dot W^j)_j$ is a $d$-dimensional white noise and $(\sigma_{i,j})_{1\le i,j\le d}$ is a non-singular deterministic matrix. 
Assume: (i) $O=(0,b)$ 
if $d=1$; (ii)  $O=B_1(0)$ if $d=2,3$.
\begin{enumerate}
\item Case $d=1,3$. From \cite[Lemmas 5.4 and 5.7]{ss-viles} we see that the random field solution to \eqref{ex-poisson}, $M(x) = (u_i(x))_{i}$,  satisfies \eqref{canon}  with $\nu=1$ and $\nu=1/2$ for $d=1$ and $d=3$, respectively. This is condition 4 of hypotheses $(H_M)$. As for conditions $1-3$ of $(H_M)$, they are established in \cite{ss-viles} (see (25) on p. 1871, Lemma 5.1 and the  proof of Theorem 5.1, respectively). Let $D>1$ if $d=1$ and $D>6$ if $d=3$.  Therefore from \eqref{hgst-example} and \eqref{hgst-example-lb}, we recover \cite[Theorems 5.10 and 5.11]{ss-viles} (with $d:=k$ there), respectively.
  
  \item Case $d=2$.  
 Let $r_0>0$ be such that $\bar B_{r_0}(0)$ is strictly contained in $O$. Claim 1 in \cite[Lemma 5.5]{ss-viles} states that, there exists a constant $C$ (depending on $r_0$ ) and for all $x,y\in B_{r_0}(0)$ with $|x-y|\le e^{-1}$, 
\beq
 \label{canon-poiss}
\left\Vert u(x) - u(y)\right\Vert_{L^2(\Omega)}\le C |x-y|\ |\log |x- y||.
\eeq
 Thus, applying Remark \ref{s2.2-r10} (with $M=u$), we deduce that  \eqref{eq2.2.17} holds with 
$q(\tau)= \tau \log\left(\tfrac{c}{\tau}\right)$. The function $q$ is of the form considered in Lemma  \ref{s2.2-l-grow-g-st} (3.) with $\nu=\delta = 1$. Furthermore, if $D>2$, $\tau\mapsto \tau \log\left(\tfrac{c}{\tau}\right)$ is an increasing function on the interval $(0,c\exp(-D/(D-2)))$. Therefore, applying Theorem \ref{2.2-tconclusive} we obtain
\beq
\label{upper-poisson-d=2}
P\left(M(K)\cap A\ne \emptyset\right)\le  C(K,\sigma_K,d,D)\ \mathcal{H}_{g_q}(A),
\eeq
for any compact set $K\subset \bar B_{r_0}(0)$, with $g_q(\tau)= \tau \log\left(\tfrac{c}{\tau}\right)$. Comparing with \cite[Theorem 5.10]{ss-viles}, we see that \eqref{upper-poisson-d=2} provides a sharper estimate.
\end{enumerate}

   \end{exa}


\section{A linear heat equation with fractional noise}
\label{s3}

We consider a non-negative definite distribution in $\mathcal{S}^\prime(\mathbb{R}^d)$ given by an absolutely continuous measure $\Lambda(dx)= f(x) dx$.
Let 
$\mu(d\xi)=\left(\mathcal{F}^{-1}f\right)(\xi) d\xi$; by the Bochner-Schwarz theorem, the measure $\mu$ is non-negative, tempered and symmetric; it is called {\em spectral measure}.  We assume that  for any non-negative measurable function $h$,
\begin{equation}
\label{eq3.3}
\int_{\mathbb{R}^d} h(\xi)\mu (d\xi)\asymp\int_{\mathbb{R}^d} h(\xi) \vert \xi\vert^{-\alpha}d\xi,\quad \text{for some}\  \alpha\in[0,d). 
\end{equation}
In an abridged form, we will write this property as 
$\mu (d\xi)\asymp \vert \xi\vert^{-\alpha}d\xi$.

Fix $\alpha\in[0,d)$, $H\in(0,1)$ and let $\{W^{H,\alpha}(t,A), t\in \left[ 0,T\right], A\in \mathcal{B}(\mathbb{R}^d) \}$ be a centred Gaussian field with covariance
\begin{equation}\label{eq3.1}
E(W^{H,\alpha}(t,A)W^{H,\alpha}(s,B))=R_H(t,s)\int_A\int_B f(z-w)dzdw,
\end{equation}
where $R_H(t,s):=\frac{1}{2}(t^{2H}+s^{2H}-\vert t-s\vert^{2H})$ is the covariance of a fractional Brownian motion with Hurst index $H$. In  this section, we restrict to the case $H\in(1/2,1)$.

If $\alpha>0$, $W^{H,\alpha}$ is called a \textit{fractional-colored} noise because it is a fractional Brownian motion in time and has a non trivial spatial covariance. 
Consider the particular case $f(x)= \delta_{\{0\}}(x)$. Then, $\mu(d\xi)= d\xi$ and \eqref{eq3.3} trivially holds with $\alpha=0$. This corresponds to the {\em fractional-white} noise, whose covariance according to \eqref{eq3.1} is 
\beqn
E(W^{H,0}(t,A)W^{H,0}(s,B))=R_H(t,s)\vert A\cap B\vert, 
\eeqn
where $\vert \cdot\vert$ denotes the Lebesgue measure.

The Riesz and the Bessel kernels are examples of functions $f$ that satisfy the above assumptions (see e.g. \cite[Ch. V]{Stein}).

Consider the linear stochastic heat equation
\beq
\frac{\partial v}{\partial t} = \Delta v+\dot{W}^{H,\alpha}, \ (t,x)\in (0 ,T]\times\mathbb{R}^d;\quad
v(0,x) = v_0(x),\  x\in \mathbb{R}^d.
\label{s3.spde}
\eeq
The random field solution to this equation is the Gaussian stochastic process 
\beq
\label{s3.1.1}
v(t,x) = I_0(t,x) + u(t,x),\quad (t,x)\in (0 ,T]\times\mathbb{R}^d,
\eeq 
where
\begin{equation}
\label{eq3.5}
I_0(t,x) = \int_{\R^d} G(t,x-y) v_0(y) dy,\quad
u(t,x) = \int_0^t\int_{\mathbb{R}^d} G(t-s,x-y)W^{H,\alpha}(ds,dy),
\end{equation}
with 
$G(t,x) = \frac{1}{(4\pi t)^{d/2}}\exp\left(-\frac{\vert x\vert^2}{4t}\right)\mathbbm{1}_{\{t\geq 0\}}$.

Suppose that the function $G(t,x-\cdot) v_0(\cdot)$ belongs to  $L^1(\R^d)$, to ensure that $x\mapsto I_0(t,x)$ is well defined for all $t\in[0,T]$. Furthermore, assume
\begin{equation}\label{eq3.6}
\int_{\mathbb{R}^d}\frac{\mu(d\xi)}{(1+\vert\xi\vert^2)^{2H}} <\infty. 
\end{equation}
Owing to \cite[Theorem 2.5]{tud}  (see also \cite[Sec. 2]{bal:tud}), this is a necessary and sufficient condition for $(u(t,x))$ given in \eqref{eq3.5} to define a $L^2(\Omega)$ random field, 
and in this case, 
$
\sup_{(t,x)\in[0,T]\times\mathbb{R}^d} E(u(t,x)^2)<\infty$.
Assuming \eqref{eq3.3}, we can check that  $(\ref{eq3.6})$ holds if and only if  $0<d-\alpha<4H$. In the remaining of the section, we will assume this constraint.
\medskip

Throughout this section, we will make use the following expression for the variance of $u(t,x), (t,x)\in(0,\infty)\times\mathbb{R}^d$: 
\begin{align}
\label{eq3.8}
\sigma^2_{t,x}:= E\left(\vert u(t, x)\vert^2\right)
 &= \alpha_H\int_0^t d\tau\int_0^t d\sigma \vert \tau-\sigma\vert^{2H-2}\notag\\
&\qquad\times \int_{\mathbb{R}^d}dz  \int_{\mathbb{R}^d}dw f(z-w)G(t-\tau,x-z)G(t-\sigma,x-w)
\notag\\
& = \frac{\alpha_H}{(2\pi)^d} \int_0^t d\tau\int_0^t d\sigma\  \vert \tau-\sigma\vert^{2H-2}\int_{\mathbb{R}^d}\mu(d\xi) 
 e^{-2(\tau+\sigma)|\xi|^2}.
\end{align}
where $\alpha_H=H(2H-1)$. The first equality can by found in \cite[Sec 2.5.1]{tud}, while the second one follows from Parseval's identity, since $\mathcal{F}(G(t,\cdot))(\xi)=e^{-t\vert \xi\vert^2}$.
From the second equality in \eqref{eq3.8}, we see that $u$ is stationary in $x$ ($\sigma^2_{tx}$ does not depend on $x$). 

For its further use, we prove some properties relative to $\sigma^2_{t,x}$.

\begin{lem}
\label{l3.1}
\begin{enumerate}
\item For any $0<t_0<T$, there exist $0<c<C<\infty$ such that, for any $(t,x)\in[t_0,T]\times \R^d$,
$c\le \sigma^2_{t,x} \le C$.
\item For any $(t,x)\in(0,\infty)\times \R^d$, the mapping 
$t\mapsto \sigma^2_{t,x}$ is differentiable.
\end{enumerate}
\end{lem} 
\proof

Use \eqref{eq3.3} in the last expression of the array \eqref{eq3.8} and then, the change of variable $\xi\mapsto (\tau+\sigma)^{\frac{1}{2}} \xi$ along with \eqref{eq3.6}. Applying the change of variables , $\tau\mapsto\frac{\tau}{t}$, $\sigma\mapsto\frac{\sigma}{t}$, we see that $\sigma_{t,x}^2$ is bounded from below (respectively, from above) by 
\beq
\label{ul-sigma-bis}
c_{\alpha,d,H} \int_0^t d\tau\int_0^t d\sigma\ \frac{\vert \tau-\sigma\vert^{2H-2}}{(\tau+\sigma)^{(d-\alpha)/2}}= t^{2H-(d-\alpha)/2}c_{\alpha,d,H}\int_0^1 d\tau\int_0^1 d\sigma\ \frac{\vert \tau-\sigma\vert^{2H-2}}{(\tau+\sigma)^{(d-\alpha)/2}}.
\eeq

Let $C_{\alpha,d,H} = \int_0^1 d\tau \int_0^1 d\sigma\frac{|\tau-\sigma|^{2H-2}}{(\tau+\sigma)^{\frac{d-\alpha}{2}}}$ and observe that, since $4H-(d-\alpha)>0$,  $C_{\alpha,d,H} <\infty$. From the above computations, we deduce that the lower inequality (respectively, the upper inequality) in the first claim holds with $c\le t^{2H-\frac{(d-\alpha)}{2}}_0 C_{\alpha,d,H}$
 (respectively, with $C\ge T^{2H-\frac{(d-\alpha)}{2}} C_{\alpha,d,H}$).

Claim 2. follows from the expression \eqref{eq3.8}.
\qed

\subsection{Equivalence for the canonical metric}\label{s3.1}
\medskip
The canonical pseudo-distance associated with the process $u$ is defined by 
\begin{equation}\label{canonical-n}
\textgoth{d}((t,x),(s,y)) =\left\Vert u(t,x)-u(s,y)\right\Vert_{L^2(\Omega)}.
\end{equation}
The goal is to prove Theorem \ref{t3.1}, which gives an equivalent pseudo-distance for $\textgoth{d}$.

We start by recalling some related results. According to 
\cite[Theorems 2.2 and 2.6]{tud}, there exist positive constants $c_1,c_2$, which depend on $\alpha$, $d$, $H$, and $T$, such that for all $t,s\in [0,T]$ and  $x\in\mathbb{R}^d$,
\begin{align}
\label{eq3.1.1}
\left\Vert u(t,x)-u(s,x)\right\Vert_{L^2(\Omega)}^2\asymp \vert t-s\vert^{2H-\frac{d-\alpha}{2}}.
\end{align}
If $\alpha\in(0,d)$, according to \cite[Theorem 4]{tud:xia}, for any fixed $t_0\in (0, T]$,  there exist  positive constants $c_3$, $c_4$ such that for any  $t\in [t_0, T]$, $x,y\in[-M,M]^d$,
\beq
\label{eq3.1.1-space}
 \left\Vert u(t,x)-u(t,y)\right\Vert_{L^2(\Omega)}^2
\asymp\left(\log\frac{1}{\vert x-y\vert}\right)^\beta\vert x-y\vert^{2\wedge (4H-(d-\alpha))}.
\eeq
where $\beta = 1$, if $4H-(d-\alpha)\geq 2$, and $\beta=0$, otherwise.

Let $W^\alpha$ be a centered Gaussian process with covariance 
\beqn
E(W^{\alpha}(t,A)W^{\alpha}(s,B))=(t \wedge s)\int_A\int_B f(z-w)\ dzdw.
\eeqn
The stochastic integral in \eqref{eq3.5} can be written as an integral with respect to $W^\alpha$(see e.g. \cite[(2.31)]{tud}):
\begin{align}
\label{representation}
&\int_0^t\int_{\R^d} G(t-r,x-z)\ W^{\alpha,H}(dr,dz)\notag\\
&\qquad \quad= \int_{\R_+} \int_{\R^d}\left(\int_{\R} d\tau\ G(t-\tau,x-z) (\tau-r)_+^{H-\frac{3}{2}}\right) W^\alpha(dr,dz).
\end{align}
Using this property, we generalize the lower bound in \eqref{eq3.1.1}, as follows.
\begin{prop}
\label{p3.1}
There exists a positive constant $c_1$ which depends on $\alpha$, $d$, $H$, and $T$, such that, for all $t,s\in[0,T]$ and $x,y\in\R^d$,
\beq
\label{eq3.1.1bis}
\Vert u(t,x)-u(s,y)\Vert^2_{L^2(\Omega)}\ge c_1\vert t-s\vert^{2H-\frac{(d-\alpha)}{2}}.
\eeq
\end{prop}
\proof
Assume, without loss of generality, that $0\le s<t\le T$. Then, from \eqref{representation}, the It\^o isometry (see \cite[Sec. 2.3.1]{tud}), and Parseval's identity, we obtain,
\begin{align}
\label{repremean}
&\Vert u(t,x)-u(s,y)\Vert^2_{L^2(\Omega)}\notag\\
&\ = E \Big(\Big\vert \int_{\R_+}  \int_{\R^d} \Big(\int_{\R} d\tau\ \left[G(t-\tau,x-z)1_{(\tau\le t)} - G(s-\tau,y-z)1_{(\tau\le s)}\right](\tau-r)_+^{H-\frac{3}{2}}\Big)\notag\\
&\qquad\qquad \times W^\alpha(dr,dz)\Big\vert^2\Big)\notag\\
&\ =\int_{\R_+} dr  \int_{\R^d} dz\int_{\R^d} dw\notag\\
&\ \left(\int_{\R} d\tau\ \left[G(t-\tau,x-z)1_{(\tau\le t)} - G(s-\tau,y-z)1_{(\tau\le s)}\right](\tau-r)_+^{H-\frac{3}{2}}\right)\notag\\
&\qquad \times \left(\int_{\R} d\tau\ \left[G(t-\tau,x-w)1_{(\tau\le t)} - G(s-\tau,y-w)1_{(\tau\le s)}\right](\tau-r)_+^{H-\frac{3}{2}}\right) f(z-w)\notag\\
&\ =(2\pi)^{-d}\int_{\R_+} dr  \int_{\R^d}\mu (d\xi)\notag\\
&\qquad \left\vert \mathcal{F}\left(\int_{\R} d\tau\ \left[G(t-\tau,x-\cdot)1_{(\tau\le t)} - G(s-\tau,y-\cdot)1_{(\tau\le s)}\right](\tau-r)_+^{H-\frac{3}{2}}\right)(\xi)\right\vert^2\notag\\
 &\ =(2\pi)^{-d}\int_{\R_+} dr  \int_{\R^d}\mu (d\xi)\left\vert\int_{\R} d\tau\ \left[e^{-i\xi x}e^{-(t-\tau)\vert\xi\vert^2}1_{(\tau\le t)} - e^{-i\xi y}e^{-(s-\tau)\vert\xi\vert^2}1_{(\tau\le s)}\right]\right.\notag\\
 &\left.\qquad \times (\tau-r)_+^{H-\frac{3}{2}}\right\vert^2.
\end{align} 
Split the domain of integration of the variable $r$ into the subdomains $[s,t]$ and $[s,t]^c$, and observe that on $[s,t]$, the term $1_{(\tau\le s)}(\tau-r)_+$ equals zero. Since the integrand is non negative, we have
\begin{align}
\label{repremeanbis}
&\Vert u(t,x)-u(s,y)\Vert^2_{L^2(\Omega)}\notag\\
&\qquad \ge(2\pi)^{-d}\ \int_s^t dr \int_{\R^d} \mu(d\xi) \left(\int_{\R} d\tau\ e^{-(t-\tau)\vert\xi\vert^2} 1_{(\tau\le t)} (\tau-r)_+^{H-\frac{3}{2}}\right)^2.
\end{align} 
Computing the integrals, we see that this is bounded below by a constant multiple of $(t-s)^{2H-\frac{d-\alpha}{2}}$, where the constant depends on $\alpha, d$ and $H$.
\qed

The next proposition extends \eqref{eq3.1.1-space} to cover the range $\alpha\in[0,d)$. The proof is the same as that of  \cite[Theorem 4]{tud:xia}, where $\alpha\in(0,d)$. For the sake of completeness, we provide the details and see that the arguments can be adapted to cover the case $\alpha=0$.


\begin{prop}\label{p3.2} Let $M>0$. There exists positive constants $c_3, c_4$, that depend on $\alpha,d, H, M$, such that for any  $t>0$, $x,y\in[-M,M]^d$,
\begin{align}
\label{eq3.1.2}
&c_3(t^{2H}\wedge 1)\left(\log\frac{2e\sqrt{d}M}{\vert x-y\vert}\right)^\beta\vert x-y\vert^{2\wedge (4H-(d-\alpha))}\notag\\
&\qquad\leq \Vert u(t,x)-u(t,y)\Vert^2_{L^2(\Omega)} \leq c_4(t^{2H}+1)\left(\log\frac{2e\sqrt{d}M}{\vert x-y\vert}\right)^\beta\vert x-y\vert^{2\wedge (4H-(d-\alpha))},
\end{align}
where $\beta=1$, if $4H-(d-\alpha)=2$, and $\beta=0$, otherwise.
\end{prop}

\proof
Similarly as in \eqref{eq3.8}, using Parseval's identity, we have
\begin{align}
\label{eq3.1.3}
&\Vert u(t,x)-u(t,y)\Vert^2_{L^2(\Omega)}\notag\\
& \quad =\frac{\alpha_H}{(2\pi)^d}  \int_0^t d\tau\int_0^t d\sigma\ \vert \tau-\sigma\vert^{2H-2} 
\int_{\mathbb{R}^d}\mu(d\xi) \ e^{-2(\tau+\sigma)|\xi|^2}\left(1- \cos[(x-y)\cdot \xi]\right).
\end{align}
According to \cite[Prop 4.3]{bal:tud2}, there exist positive constants $c_{1,H}$, $c_{2,H}$ such that
\begin{align}
\label{eq3.1.4}
c_{1,H}(t^{2H}\wedge 1)\bigg( \frac{1}{1+\vert\xi\vert^2}\bigg)^{2H}
 &\leq\int_0^t d\tau\int_0^t d\sigma \vert \tau-\sigma\vert^{2H-2} e^{-2(\tau+\sigma)|\xi|^2}\notag\\
 &\leq c_{2,H}(t^{2H}+ 1)\bigg( \frac{1}{1+\vert\xi\vert^2}\bigg)^{2H}.
\end{align}
Recall \eqref{eq3.3}. After having applied the change of variables $\xi\mapsto\frac{\eta}{\vert x-y\vert}$, from \eqref{eq3.3}, \eqref{eq3.1.3} and \eqref{eq3.1.4} we deduce
\begin{align}
\label{eq3.1.5}
&c_{1,d,H}(t^{2H}\wedge 1)\vert x-y\vert^{4H-d-\alpha} \int_{\mathbb{R}^d} d\eta \frac{\left(1-\cos\left[\left(\frac{x-y}{\vert x-y\vert}\right)\cdot\eta\right]\right)}{\vert \eta\vert^\alpha(\vert x-y\vert^2+\vert \eta\vert^2)^{2H}}
\leq\Vert u(t,x)-u(t,y)\Vert^2_{L^2(\Omega)}\notag\\
&\ \leq c_{2,d,H}(t^{2H}+ 1)\vert x-y\vert^{4H-d-\alpha} \int_{\mathbb{R}^d} d\eta \frac{\left(1-\cos\left[\left(\frac{x-y}{\vert x-y\vert}\right)\cdot\eta\right]\right)}{\vert \eta\vert^\alpha(\vert x-y\vert^2+\vert \eta\vert^2)^{2H}},
\end{align}
for all $t>0$ and $x,y\in\R^d$, with some positive and finite constants $c_{1,d,H}$, $c_{2,d,H}$. 

Next we give lower and upper bounds for the terms on the left hand side and the right hand side of \eqref{eq3.1.5}, respectively.
\smallskip

\noindent{\em Lower bounds.}  By Schwarz's inequality, $\overline{B_1(0)}\subset\left\{\eta\in\R^d: \left\vert\frac{(x-y)}{|x-y|}\cdot \eta\right\vert\le 1\right\}$. Moreover, for $|\theta|\le 1$, $1-\cos\theta\ge \frac{\theta^2}{4}$. Consequently,
\begin{align}
\label{eq3.1.6}
\mathcal{I}:= \int_{\mathbb{R}^d} d\eta \frac{\left(1-\cos\left[\left(\frac{x-y}{\vert x-y\vert}\right)\cdot\eta\right]\right)}{\vert \eta\vert^\alpha(\vert x-y\vert^2+\vert \eta\vert^2)^{2H}}\ge \frac{1}{4}\int_{\overline{B_1(0)}} d\eta \frac{\left(\frac{x-y}{\vert x-y\vert}\cdot\eta\right)^2}{\vert \eta\vert^\alpha(\vert x-y\vert^2+\vert \eta\vert^2)^{2H}}.
\end{align}
Shrink the ball $\overline{B_1(0)}$ to the spherical sector defined by the constraint $\varphi\in[0,\pi/4]$ on the angle. Then, pass to spherical coordinates and, without loss of generality, suppose that $(x-y)/|x-y|$ is the unit vector $(1,0,\ldots,0)$ in $\R^d$.
Since $\frac{x-y}{\vert x-y\vert}\cdot\eta = |\eta| \cos\varphi$, where $\varphi\in[0,\pi/4]$ is the angle between $(x-y)/|x-y|$ and $\eta$, we obtain,
\beqn
\mathcal{I}\ge C \int_0^1 d\rho \frac{\rho^{d-\alpha+1}}{\left(|x-y|^2 + \rho^2\right)^{2H}}.
\eeqn
We estimate this integral by distinguishing three cases. 
\smallskip

\noindent{\em Case 1. $0<4H-(d-\alpha)<2$.} Since $|x-y|^2 + \rho^2 \le 4dM^2 + 1$, 
\beqn
\int_0^1 d\rho \frac{\rho^{d-\alpha+1}}{(\vert x-y\vert^2+\rho^2)^{2H}}\geq\int_0^{1}d\rho \frac{\rho^{d-\alpha+1}}{(4dM^2+1)^{2H}}=\frac{1}{(d-\alpha+2)(4dM^2+1)^{2H}}.
\eeqn
\noindent{\em Case 2. $4H-(d-\alpha)=2$.} Because $|x-y|\le 2\sqrt d M$, we clearly have
\begin{align*}
\int_0^1d\rho\  \frac{\rho^{d-\alpha+1}}{(\vert x-y\vert^2+\rho^2)^{2H}}&\geq c_{\alpha,d,H,M}\int_{\frac{\vert x-y\vert}{2e\sqrt{d}M}}^1 d\rho\ \rho^{d-\alpha-4H+1}\\
&=c_{\alpha,d,H,M}\log\left(\frac{2e\sqrt{d}M}{\vert x-y\vert}\right).
\end{align*}
\noindent{\em Case 3. $4H-(d-\alpha)>2$.} Using a similar argument as for case 2,
\begin{align*}
\int_0^1d\rho \frac{\rho^{d-\alpha+1}}{(\vert x-y\vert^2+\rho^2)^{2H}}&\geq c_{\alpha,d, H, M}\int_{\frac{\vert x-y\vert}{2e\sqrt{d}M}}^1 d\rho \rho^{d-\alpha+1-4H}\\
&=c_{\alpha,d,H,M}\vert x-y\vert^{d-\alpha-4H+2}.
\end{align*}
\smallskip

\noindent{\em Upper bounds.} Apply the inequality $1-\cos(\theta)\leq 2\wedge\theta^2$ and then, use spherical coordinates to see that the integral $\mathcal{I}$ defined in \eqref{eq3.1.6} satisfies
\beq
\label{eq3.1.7}
\mathcal{I} \le\int_{\mathbb{R}^d} d\eta \frac{(2\wedge \vert\eta\vert)^2)}{\vert\eta\vert^\alpha(\vert x-y\vert^2+\vert \eta\vert^2)^{2H}}=c_d\int_0^{\infty} d\rho  \frac{(1\wedge \rho^2)\rho^{d-\alpha-1}}{(\vert x-y\vert^2+ \rho^2)^{2H}}:= c_d\ \mathcal{J}.
\eeq
We estimate $\mathcal{J}$ by considering three cases, as we did for the lower bounds.
\smallskip

\noindent{\em Case 1. $0<4H-(d-\alpha)<2$.} Since $\vert x-y\vert^2+ \rho^2\ge \rho^2$, we have
\beqn
J \leq \int_0^1 d\rho\  \rho^{d-\alpha-4H+1}+ \int_1^{\infty} d\rho\  \rho^{d-\alpha-4H-1}=c_{\alpha,d,H}.
\eeqn
\smallskip

\noindent{\em Case 2. $4H-(d-\alpha)=2$.} Splitting the domain of integration of $\mathcal{J}$, we obtain
\begin{align*}
\mathcal{J}&\leq \int_0^{\vert x-y\vert} d\rho \frac{ \rho^{d-\alpha+1}}{\vert x-y\vert^{4H}}+\int_{\vert x-y\vert}^{2e\sqrt{d}M} d\rho\  \rho^{d-\alpha-4H+1}+\int_{2e\sqrt{d}M} ^{\infty} d\rho\  \rho^{d-\alpha-4H-1}\notag\\
& =\frac{1}{(d-\alpha+2)}+\log\left(\frac{2e\sqrt{d}M}{\vert x-y\vert}\right)+\frac{(2e\sqrt{d}M)^{2}}{2}
\leq c_{\alpha,d,H,M}\log\left(\frac{2e\sqrt{d}M}{\vert x-y\vert}\right).
\end{align*}
\smallskip

\noindent{\em Case 3. $4H-(d-\alpha)>2$.} Using the inequalities $1/(|x-y|^2+\rho^2)\le 1/(|x-y|^2)$ and 
$1/(|x-y|^2+\rho^2)\le 1/\rho^2$, on $\{0\le \rho\le |x-y|\}$ and $\{|x-y|<\rho<\infty\}$, respectively, we have
\beqn
\mathcal{J}\leq \vert x-y\vert^{-4H} \int_0^{\vert x-y\vert} d\rho\  \rho^{d-\alpha+1}+ \int_{\vert x-y\vert}^{\infty} d\rho\  \rho^{d-\alpha-4H+1}
=c_{\alpha,d,H}\vert x-y\vert^{d-\alpha-4H+2}.
\eeqn

From \eqref{eq3.1.5}, and  using the lower and upper bounds obtained before, we deduce \eqref{eq3.1.2}.
 \qed
 
 We end this section by proving the equivalence for the canonical pseudo-distance \eqref{canonical-n}. It is a consequence of \eqref{eq3.1.1} and Proposition \ref{p3.2}.

\begin{teo}\label{t3.1} Fix $M>0$ and $t_0\in (0, T]$. There exists positive constants $c_5, c_6$  depending on $\alpha, d, t_0, H, M, T$ such that for any  $t,s\in [t_0, T]$ and $x,y\in[-M,M]^d$,
\begin{align}
\label{ul}
\Vert u(t,x)-u(s,y)\Vert^2_{L^2(\Omega)}\asymp
\vert t-s\vert^{2H-\frac{d-\alpha}{2}}+\left(\log\frac{2e\sqrt{d}M}{\vert x-y\vert}\right)^\beta\vert x-y\vert^{2\wedge (4H-(d-\alpha))},
\end{align}
where $\beta = 1$, if $4H-(d-\alpha)=2$, and $\beta=0$, otherwise.
 
The upper bound holds for any $t,s\in[0,T]$.
\end{teo}

\proof
The estimate from above is a consequence of the upper bounds in (\ref{eq3.1.1}) and (\ref{eq3.1.2}), which hold  for any $t,s\in[0,T]$.

We prove the estimates from below by distinguishing two cases.

\noindent{\em Case 1.\ $\vert t-s\vert ^{2H-\frac{d-\alpha}{2}}<\frac{c_3(t_o^{2H}\wedge 1)}{4c_2}\left(\log\frac{2e\sqrt{d}M}{\vert x-y\vert}\right)^\beta\vert x-y\vert^{2\wedge (4H-(d+\alpha))}$}. \
Applying the triangle inequality and then, using the lower bound in \eqref{eq3.1.2}  and the upper bound in \eqref{eq3.1.1}, we obtain
\begin{align*}
&\Vert u(t,x)-u(s,y)\Vert^2_{L^2(\Omega)}\geq \frac{1}{2} \Vert u(t,x)-u(t,y)\Vert^2_{L^2(\Omega)}-\Vert u(t,y)-u(s,y)\Vert^2_{L^2(\Omega)}\\
&\qquad\geq \frac{c_3(t_0^{2H}\wedge 1)}{2}\left(\log\frac{2e\sqrt{d}M}{\vert x-y\vert}\right)^\beta\vert x-y\vert^{2\wedge (4H-(d-\alpha))}-c_2\vert t-s\vert ^{2H-\frac{d-\alpha}{2}}\\
&\qquad\geq \frac{c_3(t_0^{2H}\wedge 1)}{8}\left(\log\frac{2e\sqrt{d}M}{\vert x-y\vert}\right)^\beta\vert x-y\vert^{2\wedge (4H-(d-\alpha))}+\frac{c_2}{2}\vert t-s\vert ^{2H-\frac{d-\alpha}{2}}.
\end{align*}

\noindent{\em Case 2.\ $\vert t-s\vert ^{2H-\frac{d-\alpha}{2}}\geq\frac{c_3(t_o^{2H}\wedge 1)}{4c_2}\left(\log\frac{2e\sqrt{d}M}{\vert x-y\vert}\right)^\beta\vert x-y\vert^{2\wedge (4H-(d-\alpha))}$}. \
By Proposition \ref{p3.1},
\begin{align*}
&\Vert u(t,x)-u(s,y)\Vert^2_{L^2(\Omega)}\geq c_1\vert t-s\vert ^{2H-\frac{d-\alpha}{2}}\\
&\qquad\geq \frac{c_1}{2}\vert t-s\vert ^{2H-\frac{d-\alpha}{2}}+\frac{c_3(t_o^{2H}\wedge 1)}{8c_2}\left(\log\frac{2e\sqrt{d}M}{\vert x-y\vert}\right)^\beta\vert x-y\vert^{2\wedge (4H-(d-\alpha))}.
\end{align*}
The proof is complete.
\qed

\begin{rem}
\label{s3.1-r30}
Assume that $v_0\in \mathcal{C}^\zeta(\R^d)$, for some $\zeta\in (0,1]$. Then the function 
\beqn
[0,T] \times \R^d\ni (t,x) \longrightarrow I_0(t,x) = \int_{\R^d} G(t,x-y) v_0(y) dy,
\eeqn
is globally H\"older continuous, jointly is $(t,x)$, with exponents $(\zeta/2, \zeta)$ (see e.g. \cite{d-ssBook}).

Furthermore, the upper bound estimate in \eqref{ul} and the classical Kolmogorov's continuity criterion ensures the existence of a version of the process $(v(t,x))$ with continuous (and even H\"older continuous) sample paths, jointly in $(t,x)$.
\end{rem}

We end this section giving some properties of the covariance function of the process $(u(t,x))$ that will be used in Section \ref{s3.2}.

\begin{lem}
\label{l3.2}
Fix $M>0$ and $t_0\in(0,T]$. 
\begin{enumerate}
\item There exists $\eta>0$ and $C>0$, depending on $\alpha,d, t_0, H, M, T$, such that, for all $s,t\in[t_0, T]$ and $x,y\in[-M,M]^d$,
\beq
\label{eq3.2.2}
\vert\sigma^2_{t,x}-\sigma^2_{s,y}\vert\leq C\ \Vert u(t,x)- u(s,y)\Vert_{L^2(\Omega)}^{1+\eta}.
\eeq
\item For any $(t,x), (s,y)\in[t_0,T]\times \R^d$ such that $(t,x)\ne (s,y)$, 
\beqn
\rho_{(t,x),(s,y)}<1. 
\eeqn
\end{enumerate}
\end{lem}
\proof
1. Assume, without loss of generality, that $0<s\le t$. For all $x,y\in\R^d$, from \eqref{eq3.8}  and  similarly as in \eqref{ul-sigma-bis}, we deduce
\begin{align}\label{dif-sigmas}
&\left(\frac{\alpha_H}{(2\pi)^d}\right)^{-1}\vert\sigma^2_{t,x}-\sigma^2_{s,y}\vert = \left(\frac{\alpha_H}{(2\pi)^d}\right)^{-1}\left(\sigma^2_{t,x}-\sigma^2_{s,y}\right)\notag\\
&\ = \int_{\mathbb{R}^d}\mu(d\xi)e^{-2(\tau+\sigma)|\xi|^2}\bigg( \int_0^t d\tau\int_0^t d\sigma\  \vert \tau-\sigma\vert^{2H-2}-  \int_0^s d\tau\int_0^s d\sigma\  \vert\tau-\sigma\vert^{2H-2}\bigg)\notag\\
& \ \leq c_{\alpha,d,H}\left(\int_s^t d\tau\int_s^t d\sigma\ \frac{\vert\tau-\sigma\vert^{2H-2}}{(\tau+\sigma)^{\frac{d-\alpha}{2}}}+2\int_0^s d\tau \int_s^t d\sigma\ \frac{\vert\tau-\sigma\vert^{2H-2}}{(\tau+\sigma)^{\frac{d-\alpha}{2}}}\right).
\end{align}

Apply polar coordinates $(\tau,\sigma)\mapsto(\rho\cos\theta,\rho\sin\theta)$ and then, the mean value theorem, to see that
\begin{align*}
&\int_s^t d\tau\int_s^t d\sigma\ \frac{\vert \tau-\sigma\vert^{2H-2}}{(\tau+\sigma)^{\frac{d-\alpha}{2}}}\\
&\quad \leq 
\int_{\sqrt{2}s}^{\sqrt{2}t} d\rho\ \rho^{2H-\frac{d-\alpha}{2}-1}
\left(\int_{0}^{\frac{\pi}{2}}d\theta
\frac{\vert \cos\theta-\sin\theta\vert^{2H-2}}{(\cos\theta+\sin\theta)^{\frac{d-\alpha}{2}}}\right)
\notag\\
&\quad \leq \frac{2^{H-\frac{(d-\alpha)}{4}}T^{2H-\frac{(d-\alpha)}{2}-1}(t-s)}{\left(2H-\frac{(d-\alpha)}{2}\right)^2} \int_{0}^{\frac{\pi}{2}}d\theta\ \frac{\vert \cos\theta-\sin\theta\vert^{2H-2}}{(\cos\theta+\sin\theta)^{\frac{d-\alpha}{2}}}
\le C(\alpha,d,H,T)(t-s).
\end{align*}
Since $0<2H-\frac{(d-\alpha)}{2}<2$, we have $\eta_1:= \left(H-\frac{(d-\alpha)}{4}\right)^{-1}-1>0$, and we deduce,
\beq
\label{eq3.2.3}
\int_s^t d\tau\int_s^t d\sigma\ \frac{\vert \tau-\sigma\vert^{2H-2}}{(\tau+v)^{\frac{d-\alpha}{2}}}\le C(H,d,T)(t-s)^{\frac{4H-(d-\alpha)}{4}(1+\eta_1)}.
\eeq
As for the second integral on the last line of \eqref{dif-sigmas}, we have
\begin{equation}
\label{eq3.2.4}
\int_0^s d\tau\int_s^t d\sigma\ \frac{\vert\tau-\sigma\vert^{2H-2}}{(\tau+\sigma)^{\frac{d-\alpha}{2}}}\leq \int_0^s d\tau\int_s^t d\sigma\ (\sigma-\tau)^{2H-\frac{(d-\alpha)}{2}-2},
\end{equation}
because $\tau\le \sigma$ implies $\tau+\sigma\ge \sigma-\tau$.

Our next goal is to obtain estimates from above on the right-hand side of \eqref{eq3.2.4} in terms of powers of $(t-s)$. For this, we consider three cases.
\smallskip

\noindent{\em Case 1. $0<4H-(d-\alpha)<2$}.
\begin{align}
\label{eq3.2.5}
 \int_0^s d\tau\int_s^t d\sigma\ &(\sigma-\tau)^{2H-\frac{(d-\alpha)}{2}-2}=\frac{s^{2H-\frac{(d-\alpha)}{2}}+(t-s)^{2H-\frac{(d-\alpha)}{2}}-t^{2H-\frac{(d-\alpha)}{2}}}{(2H-\frac{(d-\alpha)}{2})(1+\frac{(d-\alpha)}{2}-2H)}\notag\\
&\leq \frac{(t-s)^{2H-\frac{(d-\alpha)}{2}}}{(2H-\frac{(d-\alpha)}{2})(1+\frac{(d-\alpha)}{2}-2H)}=C(\alpha,d,H)(t-s)^{2H-\frac{d-\alpha}{2}}\notag\\
&= C(\alpha,d,H)(t-s)^{\frac{4H-(d-\alpha)}{4}(1+\eta_2)},
\end{align}
with $\eta_2=1$
\smallskip

\noindent{\em Case 2.  $0<4H-(d-\alpha)=2$}. 
\begin{align*} 
\int_0^s d\tau\int_s^t d\sigma\ (\sigma-\tau)^{-1}&=t\log(t)-s\log(s)+(t-s)\log\left((t-s)^{-1}\right)\\
&\leq 2[(t\log t-s\log s)\vee ((t-s)\log((t-s)^{-1}))]\\
&\leq 2(t-s)[( \log T+1) \vee \log((t-s)^{-1})],
\end{align*}
where in the last inequality we have applied the mean value theorem.
This yields, for any $\gamma\in (0,1)$,
\begin{align}
\label{eq3.2.6}
\int_0^s d\tau\int_s^t d\sigma\ (\sigma-\tau)^{-1}&\leq 2( \vert\log T\vert+2)([(t-s)^\gamma \vee (t-s)]\le C(T)\ (t-s)^\gamma\notag\\
&=C(T)(t-s)^{\frac{4H-(d-\alpha)}{4}(1+\eta_3)},
\end{align}
with $\eta_3= 2\gamma-1$
\smallskip

\noindent{\em Case 3.  $2<4H-(d-\alpha)<4$}. 
\begin{align}
 \label{eq3.2.7}
 \int_0^s d\tau\int_s^t d\sigma\ &(\sigma-\tau)^{2H-2-\frac{(d-\alpha)}{2}}=\frac{t^{2H-\frac{(d-\alpha)}{2}}-s^{2H-\frac{(d-\alpha)}{2}}-(t-s)^{2H-\frac{(d-\alpha)}{2}}}{\left(2H-\frac{(d-\alpha)}{2}\right)\left(2H-1-\frac{(d-\alpha)}{2}\right)}\notag\\
 &\le \frac{t^{2H-\frac{(d-\alpha)}{2}}-s^{2H-\frac{(d-\alpha)}{2}}}{\left(2H-\frac{(d-\alpha)}{2}\right)\left(2H-1-\frac{(d-\alpha)}{2}\right)}
 \le \frac{T^{2H-1-\frac{(d-\alpha)}{2}}}{2H-1-\frac{(d-\alpha)}{2}}(t-s)\notag\\
 &\le C(\alpha,d,H,T)(t-s)^{\frac{4H-(d-\alpha)}{4}(1+\eta_4)},
 \end{align}
 with $\eta_4=\eta_1=\left(H-\frac{(d-\alpha)}{4}\right)^{-1}-1$.
 
Set $\eta=\min(\eta_i, i=1,2,3)$. Appealing to Theorem \ref{t3.1}, and using \eqref{dif-sigmas}, \eqref{eq3.2.3}, \eqref{eq3.2.5}, \eqref{eq3.2.6} and \eqref{eq3.2.7}, we obtain
 \begin{align*}
 \vert\sigma_{t,x}^2 - \sigma_{s,y}^2\vert &\le C(\alpha,d,H,T)(t-s)^{\left(H-\frac{d-\alpha}{4}\right)(1+\eta)}
\le c_5^{-1}C(\alpha,d,H,T) \Vert u(t,x)-u(s,y)\Vert^{1+\eta}_{L^2(\Omega)}, 
 \end{align*}
 with $c_5$ as in \eqref{ul}. The proof of Claim 1. is complete.
 \medskip

Next, we prove Claim 2 of the Lemma.  Assume that $\rho_{(t,x),(s,y)}=1$ and hence, that there exists $\lambda\in\mathbb{R}$ such that 
\begin{equation}
\label{eq3.2.8} 
\Vert u(t,x)-\lambda u(s,y)\Vert_{L^2(\Omega)}=0.
\end{equation}
We will see that this assumption leads to a contradiction.

\noindent{\em Case 1. $s<t$.}\  Apply \eqref{repremean} with $u(s,y)$ replaced by $\lambda u(s,y)$ to obtain
\begin{align*}
\Vert u(t,x)-\lambda u(s,y)&\Vert_{L^2(\Omega)}^2
=(2\pi)^{-d}\int_{\R_+} dr  \int_{\R^d}\mu (d\xi)\\
&\ \times \left\vert\int_{\R} d\tau\ \left[e^{-2(t-\tau)\vert\xi\vert^2}1_{(\tau\le t)} - \lambda e^{-2(s-\tau)\vert\xi\vert^2}1_{(\tau\le s)}\right](\tau-r)_+^{H-\frac{3}{2}}\right\vert^2.
\end{align*}
As in \eqref{repremeanbis}, this is bounded from below by a constant multiple of
\beqn
\int_{\R^d} \mu(d\xi)\int_s^t dr  \left(\int_r^t d\tau\ e^{-2(t-\tau)\vert\xi\vert^2} (\tau-r)^{H-\frac{3}{2}}\right)^2.
\eeqn
A direct computation shows that $\int_s^t dr \left(\int_r^t d\tau\ e^{-2(t-\tau)\vert\xi\vert^2} (\tau-r)^{H-\frac{3}{2}}\right)^2\ne 0$. Since we are assuming \eqref{eq3.2.8}, we reach a contradiction.

We notice that, in the case under consideration, the arguments hold for any $(t,x), (s,y)\in[0,\infty) \times \R^d$.
\medskip

\noindent{\em Case 2. $s=t\in[t_0,T]$, $x\ne y$.}\  Apply \eqref{eq3.1.3} with $u(t,y)$ replaced by $\lambda u(t,y)$ to see that
\begin{align*}
\Vert u(t,x)-\lambda u(s,y)\Vert_{L^2(\Omega)}^2
&= \frac{\alpha_H}{(2\pi)^d} \int_0^t d\tau\int_0^t d\sigma\ \vert \tau-\sigma\vert^{2H-2}\\ 
&\qquad\times \int_{\mathbb{R}^d}\mu(d\xi) 
\ e^{-2(\tau+\sigma)|\xi|^2}\left(1+\lambda^2- 2\lambda\cos[(x-y)\cdot \xi]\right).
\end{align*}
Using the lower bound estimates in \eqref{eq3.3} and \eqref{eq3.1.4}, we deduce 
\begin{align*}
\Vert u(t,x)-\lambda u(s,y)\Vert_{L^2(\Omega)}^2 &\ge C(\alpha,d,t_0,H) \int_{\R^d} \left(1+\lambda^2- 2\lambda\cos[(x-y)\cdot \xi]\right)\\
&\qquad\times\frac{\vert\xi\vert^{-\alpha}}{(1+|\xi|^2)^{2H}}\ d\xi.
\end{align*}
By assumption, the integral on the right-hand side must be zero. However, this integral is bounded from below by the integral on the spherical sector of the ball $B_1(0)$ where $2\cos[(x-y)\cdot \xi] \le \frac{1+\lambda^2}{2}$. Consequently, 
\beqn
0=  \Vert u(t, x) - \lambda\ u(t,y)\Vert^2_{L^2(\Omega)}\ge C(\alpha,d,t_0,H) \frac{1+\lambda^2}{2} \int_0^1\frac{r^{d-\alpha-1}}{(1+r^2)^{2H}} dr.
\eeqn
Since $\int_0^1\frac{r^{d-\alpha-1}}{(1+r^2)^{2H}} dr>0$, this is a contradiction. This ends the proof of Claim 2.
\qed

\subsection{Hitting probabilities}
\label{s3.2} 

Consider the random field 
$U=\{U(t,x)=(U_1(t,x),...,U_D(t,x)),\ (t,x)\in[0,T]\times\R^{d}\}$,
where the components are independent copies of the random variable $v(t,x)$ defined in \eqref{s3.1.1}. The process $U$ is  the random field solution to the system of SPDEs 
\beq
\label{eq3.2.1}
\begin{cases}
\frac{\partial U_j}{\partial t}(t,x) = \Delta U_j(t,x) +\dot{W}_j^{H,\alpha}, & (t,x)\in(0 ,T]\times\mathbb{R}^d, \notag\\
U_j(0,x) = v_0(x), & x\in \mathbb{R}^d,
\end{cases}
\eeq
$j=1,\ldots,D$, where $(W_j^{H,\alpha}, j=1,\ldots, D)$ are independent copies of the fractional-colored noise $W^{H,\alpha}$ introduced at the beginning of Section \ref{s3}. We will write $U_j(t,x) = I_{0}(t,x)+ u_j(t,x)$.
In the sequel, we assume that $v_0$ is such that the function $(t,x)\mapsto I_0(t,x)$ is continuous (see Remark \ref{s3.1-r30} for sufficient conditions).

Throughout this section, we will consider the compact sets $I=[t_0,T]$ and $J=[-M,M]^d$, with $t_0\in (0,T]$, $M>0$, and the gauge functions defined in $\R_+$,
\beq
\label{defqus}
q_1(\tau) = \tau^{H-\frac{d-\alpha}{4}},\qquad
q_2(\tau) =  
\begin{cases}
\tau^{1\wedge \left(2H-\frac{d-\alpha}{2}\right)}, & \text{if}\ 4H-(d-\alpha) \ne 2,\\
\tau\left(\log\frac{2e\sqrt d M}{\tau}\right)^{\frac{1}{2}}, & \text{if}\ 4H-(d-\alpha) =2.
\end{cases}
\eeq 
\medskip

If $4H-(d-\alpha) \ne 2$, the functions $q_1$ and $q_2$ belong to the class of examples considered in Lemma \ref{s2.2-l-grow-g-st} (1.), with $\nu_1:=H-\frac{d-\alpha}{4}$, $\nu_2:= 1\wedge \left(2H-\frac{d-\alpha}{2}\right)$. If $D>\left(\frac{1}{H-\frac{d-\alpha}{4}}+\frac{d}{1\wedge\left(2H-\frac{d-\alpha}{2}\right)}\right)$, the function 
\beq
\label{1000}
\bar g_q(\tau) = \tau^{D-\left(\frac{1}{H-\frac{d-\alpha}{4}}+\frac{d}{1\wedge\left(2H-\frac{d-\alpha}{2}\right)}\right)}
\eeq
(see \eqref{gstbis})  is strictly increasing.

Furthermore, we prove in Lemma \ref{s2.1-l2} (1.) that the function $\bar v_q(\tau)$ defined in \eqref{2.1.51} satisfies the condition \eqref{2.1.52} with $\bar g_q$ given in \eqref{1000}.
\smallskip

If $4H-(d-\alpha) = 2$,  $q_1$ and $q_2$ belong to the class of examples considered in Lemma \ref{s2.2-l-grow-g-st} (2.) with 
$\nu_1:=H-\frac{d-\alpha}{4}$, $\nu_2 = 1$,  $\delta = \frac{1}{2}$. If $D>\frac{1}{H-\frac{d-\alpha}{4}}+d$, the function 
\beq
\label{10001}
\bar g_q(\tau) = \tau^{D-\frac{1}{H-\frac{d-\alpha}{4}}} \left(q_2^{-1}(\tau)\right)^{-d},
 \eeq
 is strictly increasing on a small interval $(0,\rho_0)$. Moreover, according to Lemma \ref{s2.1-l2} (3.), 
 this function satisfies the condition \eqref{2.1.52}, where $\bar v_q(\tau)$ is defined in \eqref{2.1.51}.
 \medskip

We now give the main theorem on hitting probabilities for the process $U$.
\begin{teo}
\label{t3.2-ub} 
Let $t_0>0$, $I=[t_0,T]$, $J=[-M,M]^d$. Suppose that the function $I\times J\ni(t,x)\mapsto I_0(t,x)$ satisfies the condition \eqref{upexpect}.
\begin{enumerate}
\item Case $4H-(d-\alpha) \ne 2$. Assume $D >\left(\frac{1}{H-\frac{d-\alpha}{4}}+\frac{d}{1\wedge\left(2H-\frac{d-\alpha}{2}\right)}\right)$ and let $\bar g_q$ be as in \eqref{1000}.
\begin{enumerate}
\item There exists a constant
 $C:=C(I,J, D, d)$ such that for any Borel set $A\subset {\R}^D$,
\beq
\label{u-l-b-aniso-1}
P(U(I  \times J) \cap A\neq \emptyset)\leq C\ \mathcal{H}_{\bar g_q}(A).
\eeq
\item Fix $N>0$ and let $A\subset B_N(0)\subset \R^D$ be a Borel set. There exists a constant $c:= c(I,J, N, D, d)$ such that 
\beq
\label{u-l-b-aniso-1-bis}
P(U(I  \times J) \cap A\neq \emptyset)\ge c\ \text{Cap}_{(\bar g_q)^{-1}}(A).
\eeq
\end{enumerate}
\item Case $4H-(d-\alpha) =2$. Assume $D >\frac{1}{H-\frac{d-\alpha}{4}}+d$ and let $\bar g_q$ be  as in \eqref{10001}.
\begin{enumerate}
\item There exist a constant
 $C:=C(I,J,D,d)$ such that for any Borel set $A\subset{\R}^D$,
\beq
\label{u-l-b-aniso-2}
P(U(I  \times J) \cap A\neq \emptyset)\leq C\ \mathcal{H}_{\bar g_q}(A).
\eeq
\item Fix $N>0$ and let $A\subset B_N(0)\subset \R^D$ be a Borel set. There exists a constant $c:= c(I,J, N, D, d)$ such that 
\beq
\label{u-l-b-aniso-2-bis}
P(U(I  \times J) \cap A\neq \emptyset)\ge c\ \text{Cap}_{(\bar g_q)^{-1}}(A).
\eeq
\end{enumerate}
\end{enumerate}
\end{teo}
\proof
(i)\ {\em Upper bounds.}
The inequalities \eqref{u-l-b-aniso-1} and \eqref{u-l-b-aniso-2} are obtained applying Theorem \ref{t2.2-tconclusive-st}. Indeed, the random field $U$ is Gaussian and has i.i.d. components and has continuous
sample paths, a.s. Lemma \ref{l3.1} (1.) gives the non degeneracy condition $\sigma^2_{I,J}>0$ on the variances. Furthermore from the discussion at the begining of this section, we see that the hypotheses on the gauge functions and the corresponding $\bar g_q$ are satisfied. 
Finally, the upper bound in \eqref{ul} implies the validity of \eqref{s2.2-200} and therefore, by Remark \ref{s2.2-r20}, that of \eqref{qoneqtwo}. 
Hence, in the two cases, $U$ satisfies the hypotheses of Theorem \ref{t2.2-tconclusive-st}.
\smallskip

Observe that, when  $4H-(d-\alpha) \ne 2$, if $D-1/(H-\frac{d-\alpha}{4})+d/(1\wedge(2H-\frac{d-\alpha}{2}))<0$, we have $\mathcal{H}_{\bar g}(A)=\infty$; thus \eqref{u-l-b-aniso-1} still holds but is not informative. 
\smallskip

(ii)\ {\em Lower bounds.} The inequalities \eqref{u-l-b-aniso-1-bis} and \eqref{u-l-b-aniso-2-bis} are obtained applying Theorem \ref{s2.1-t2}. For this, we first check that the process $U$  satisfies the hypotheses $(H_U)$ of Section \ref{s2.1}. Indeed, \eqref{2.1.1-bis} is Lemma \ref{l3.1} (1.), and conditions 2 and 3 are proved in Lemma \ref{l3.2}. Theorem \ref{t3.1}
tells us that \eqref{2.1.3-bisbis} is satisfied with $q_1$ and $q_2$ given in \eqref{defqus}.  Hence, $(H_U)$ is satisfied. The conditions required on the gauge functions and the corresponding functions $\bar g_q$ and $\bar v_q$ are proved in Lemmas \ref{s2.2-l-grow-g-st} and \ref{s2.1-l2} (3.). Thus, in the two cases, $U$ satisfies the hypotheses of Theorem \ref{s2.1-t2}.

The proof of the theorem is complete.
 \qed

\section{Auxiliary lemmas}
\label{appendix}

In the next lemmas, $q$, $q_1$, $q_2$ are gauge functions and $g_q$, $\bar g_q$ the functions defined in \eqref{gg}, \eqref{gstbis}, respectively.  For convenience we recall their respective expressions:
\beqn
g_q(\tau) = \frac{\tau^D}{\left(q^{-1}(\tau)\right)^d}, \quad  \bar g_q(\tau) = \frac{\tau^D}{\left(q_1^{-1}(\tau)\right)^{d_1}\left(q_2^{-1}(\tau)\right)^{d_2}},\quad \quad\tau\in \R_+,
\eeqn
where $\bar g_q$, stands for $\bar g_{(q_1,q_2)}$. Observe that if $q_1=q_2:=q$ then $\bar g_q= g_q$ with $d:=d_1+d_2$.
\medskip

\begin{lem}
\label{s2.2-l-grow-g-st}
Fix $\rho_0>0$. Assume that  $q_1$, $q_2$  are differentiable in $(0,\rho_0)$. Then $\bar g_q$ is strictly increasing on $(0,\rho_0)$ if and only if 
\beq
\label{barg-incr}
D>\tau\left(\frac{d_1}{q_1^{-1}(\tau) \dot q_1(q_1^{-1}(\tau))} + \frac{d_2}{q_2^{-1}(\tau) \dot q_2(q_2^{-1}(\tau))}\right),\quad \tau\in(0,\rho_0),
\eeq
or equivalently, if and only if for any $\tau\in (0, q_2^{-1}(\rho_0))$,
\beq
\label{barg-incr-bis}
D> q_2(\tau)\left(\frac{d_1}{q_1^{-1}(q_2(\tau)) \dot q_1(q_1^{-1}(q_2(\tau)))} + \frac{d_2}{\tau\dot q_2(\tau)}\right).
\eeq

When $q_1 = q_2:=q$, the condition \eqref{barg-incr} is
\beq
\label{g-increase-previous}
D>d\ \frac{\tau }{q^{-1}(\tau) \dot q(q^{-1}(\tau))}, \  \  \tau\in(0,\rho_0)\ \Longleftrightarrow\ 
D> d\ \frac{q(\tau)}{\tau \dot q(\tau)}, 
  \ \  \tau\in(0,q^{-1}(\rho_0)),
\eeq
whith $d=d_1+d_2$.

For the gauge functions listed below, we have the following.
\begin{enumerate}
\item Let  $q_i(\tau) = \tau^{\nu_i}$, $\tau\ge 0$, $\nu_i>0$, $i=1,2$. Assume that $D>\tfrac{d_1}{\nu_1}+\tfrac{d_2}{\nu_2}$. Then condition  \eqref{barg-incr} holds on $\R_+$ and therefore, $\bar g_q$ is strictly increasing.
Moreover, \eqref{2.c-polar-bis} is satisfied if and only if $D>\tfrac{d_1}{\nu_1}+\tfrac{d_2}{\nu_2}$.

In particular, if $q_1(\tau)=q_2(\tau)= \tau^{\nu}$ and $d_1+d_2 =d$, the function $g_q$  is strictly increasing  on $\R_+$ whenever $D>d/\nu$. The condition  \eqref{2.c-polar} is satisfied if and only if $D>d/\nu$ holds.

\item Let $q_1(\tau) = \tau^{\nu_1}$, $q_2(\tau) = \tau^{\nu_2}\left(\log\frac{c}{\tau}\right)^\delta $, $\tau\ge 0$, $\nu_1, \nu_2, \delta>0$. Assume that $D>\tfrac{d_1}{\nu_1}+\tfrac{d_2}{\nu_2}$ and $\nu_2\ge \delta$.
Set $\eta:=(\nu_1 \nu_2D-\nu_2 d_1 - \nu_1 d_2)/(\nu_1 D-d_1)$. 
Then, on the interval $(0, c\min(e^{-1},\exp(-d/\eta))$, the condition \eqref{barg-incr-bis} holds and therefore, $\bar g_q$ is strictly increasing on this interval. The condition  \eqref{2.c-polar-bis} holds if and only if $D>\tfrac{d_1}{\nu_1}+\tfrac{d_2}{\nu_2}$.

\item Let $q(\tau) = \tau^\nu\left(\log\left(c/\tau\right)\right)^\delta$, $\tau\ge 0$, $\nu, \delta>0$. Suppose $d/D<\nu$. If $\nu-\delta <d/D<\nu$ then $g_q$ is strictly increasing on  $\tau\in\left(0, c\exp\left(-\delta/(\nu-d/D)\right)\right)$. If $d/D\le \nu-\delta$ then $g_q$ is  strictly increasing  on  $\tau\in(0, c/e)$. Furthermore, the condition \eqref{2.c-polar} holds if  and only if $D>d/\nu$. \end{enumerate}
All these examples consist of infinitely differentiable functions and, if $\nu, \nu_1, \nu_2 \in(0,1)$, the first order derivatives are decreasing on  $(0,r_0)$. For $\tau\mapsto \tau^\nu$, $r_0=\infty$, while for $\tau\mapsto 
\tau^\nu\left(\log\left(c/\tau\right)\right)^\delta$, $r_0 = c\ exp(-(1-\delta)/(1-\nu))$.

\end{lem}
\proof
Imposing the constraint $\dot{\bar{g}}(\tau)>0$ for any $\tau\in(0,\rho_0)$, yields \eqref{barg-incr}.  The equivalent form \eqref{barg-incr-bis}  is obtained by the change of variable $\tau\mapsto q_2^{-1}(\tau)$. Taking $q_1=q_2=q$, yields \eqref{g-increase-previous}.

The results on monotonicity concerning the three examples can be argued by elementary computations on the expressions  \eqref{barg-incr},  \eqref{barg-incr-bis} and \eqref{g-increase-previous}, respectively.

In the examples discussed in 1. and under the given conditions, the validity of \eqref{2.c-polar} is trivial. 
Let $q$ be as in 3. The inverse $q^{-1}$ is given by the relation 
\beq
\label{2.1.35}
q^{-1}(\tau) = c \exp\left[\frac{\delta}{\nu} W_{-1}\left(-\frac{\nu}{\delta} c^{-\frac{\nu}{\delta}}
\tau^{\frac{1}{\delta}}\right)\right],
\eeq
where $W_{-1}$ is the  real branch of the multi-valued Lambert function $W(z)$ defined for $z\in(-e^{-1},-1)$, satisfying $W(z)\le -1$. Acording to \cite[Theorem 1]{chat}, 
\beq
\label{Lambertequiv}
-1 - \sqrt{2z} - z < W_{-1}(-e^{-z-1}) < -1-\sqrt{2z} -\frac{2}{3} z, \quad z>0.
\eeq
Applying this result, we see that
\beqn
q^{-1}(\tau) \asymp c_1 \tau^{\frac{1}{\nu}} \exp\left(-\sqrt 2 \frac{\delta}{\nu}\left(-\log\left(c_2 \tau^{\frac{1}{\delta}}\right)\right)^{\frac{1}{2}}\right),
\eeqn
where $c_1, c_2$ are constants depending on $\nu, \eta$; consequently.
\beqn
g_q(\tau) \asymp \tau^{D-\frac{d}{\nu}}\exp\left[C_1\left(\log\frac{1}{ c_2 \tau^{\frac{1}{\delta}}}\right)^{\frac{1}{2}}\right].
\eeqn
Assuming $D>\frac{d}{\nu}$, the limit of the right-hand side of the above equivalence tends to zero as $\tau\downarrow 0$. Therefore, \eqref{2.c-polar} holds.

With similar arguments, one checks that in Example 2., \eqref{2.c-polar-bis} holds.


Finally, after computation of the second derivatives and the analysis of their sign, we obtain the last statement.  
\qed


In the next lemma we study properties of the functions $v_q$ and $\bar v_q$ defined in \eqref{1.3-bis} and \eqref{2.1.51}, respectively, for the particular cases of gauge functions relevant to this article.

\begin{lem}
\label{s2.1-l2}
\begin{enumerate}
\item Let $q_i(\tau) = \tau^{\nu_i}$, \ $\tau\ge 0$, with $\nu_i>0$, $i=1,2$. Let $\chi = \frac{d_1}{\nu_1}+\frac{d_2}{\nu_2}$.
Then
\beq
\label{2.1.65}
\bar v_q(\tau) = \begin{cases}
(\nu_1\nu_2(D-\chi))^{-1} \left[\tau^{-(D-\chi)} - c_{I,J}^{-(D-\chi)}\right] & \text{if}\ \chi\ne D,\\
(\nu_1\nu_2)^{-1} \log\left(\frac{c_{I,J}}{\tau}\right), & \text{if}\  \chi=D.\\
\end{cases}
\eeq
Therefore, up to multiplicative constants, $\bar v_q$ is bounded above by the Bessel-Riesz potential kernel of order $\beta:= D- \chi$.

If $\chi<D$, the function $\bar g_q$, which in this particular example is
$\bar g_q (\tau)= \tau^{D-\chi}$, satisfies the condition \eqref{2.1.52}.
\smallskip

In the particular case $q(\tau):=q_1(\tau)=q_2(\tau)=\tau^\nu$, $\tau\ge 0$, $\nu>0$, we have $\chi=\tfrac{d}{\nu}$ with $d=d_1+d_2$, and 
$\bar v_q = \nu^{-1} v_q$. Therefore,
\beq
\label{2.1.32}
v_q(\tau)=
\begin{cases}
(\nu D-d)^{-1} \left[\tau^{-(D-d/\nu)} - c_{I,J}^{-(D-d/\nu)}\right], & {\text{if}}\quad d/\nu \ne D,\\
\nu^{-1} \log \left(\frac{c_{I,J}^\nu}{\tau}\right), & {\text{if}}\quad  d/\nu =D.
\end{cases}
\eeq
Hence,  if $d/\nu < D$,  the function $g_q$ satisfies \eqref{2.1.9}.
\item Let $q(\tau)= \tau^\nu \left(\log\frac{c}{\tau}\right)^\delta$, $\tau\ge 0$, with $\nu >0$, $\delta>0$. Then,
\beq
\label{2.1.33}
v_q(\tau) \asymp 1,\  \text{if either}\  d/\nu>D\  \text{or}\  (d/\nu=D, 1-\delta D<0),
\eeq
while if either $d/\nu=D, 1-\delta D\ge 0$ or $d/\nu<D$, 
\beq
\label{2.1.34}
\lim_{\tau\downarrow 0} v_q(\tau) = \infty.
\eeq
Furthermore, the function $g_q$  satisfies \eqref{2.1.9} only if $d/\nu<D$ .
\item Let $q_1(\tau) = \tau^{\nu_1}$, $q_2(\tau) = \tau^{\nu_2}\left(\log\frac{c}{\tau}\right)^\delta$, \ $\tau\ge 0$, with $\nu_i,  \delta >0$, $i=1,2$. Then, if $D>\tfrac{d_1}{\nu_1}+\tfrac{d_2}{\nu_2}$ and $d_2\ge \nu_2$, the function $\bar g_q$, which in this case is
\beq
\label{2.1.37}
\bar g_q(\tau)=\tau^{D-\tfrac{d_1}{\nu_1}} \left(q_2^{-1}(\tau)\right)^{-d_2},
\eeq
satisfies \eqref{2.1.52}.
\end{enumerate}

\end{lem}
\proof
1. 
Computing the integral \eqref{2.1.51} for the particular choice of gauge functions $q_1$, $q_2$, we obtain \eqref{2.1.65}.  Up to multiplicative constants, this is indeed bounded by the  Bessel-Riesz potential kernel of order $\beta:= D- \chi$. 

If $\chi< D$, 
$\bar v_q(\tau)\le \left((\nu_1\nu_2( D-\chi)\ \bar g_q(\tau)\right)^{-1}$, and therefore  \eqref{2.1.52} holds. Particularizing to $q(\tau):=q_1(\tau)=q_2(\tau)=\tau^\nu$ yields \eqref{2.1.32} and its consequences.
\medskip

2.\  Properties \eqref{2.1.33} and \eqref{2.1.34} are proved using  \eqref{2.1.35} and \eqref{Lambertequiv}.

Since $v_q$ and $g_q$ are continuous functions on $(0,\infty)$, the condition \eqref{2.1.9} is equivalent to $\lim_{\tau\downarrow 0} v_q(\tau) g_q(\tau)\in(0,\infty)$.  Furthermore, because $\tau\mapsto q(\tau)$ is strictly increasing and $q(0)=0$, this is equivalent to
$
\lim_{\tau\downarrow 0} v_q(q(\tau)) g_q(q(\tau)) = l_0\in(0,\infty).
$

Consider first the cases: (i) $D<d/\nu$; (ii) $D=d/\nu$ and $1-\delta D<0$.  Since $v_q\asymp 1$ and $\lim_{\tau\downarrow 0} g_q(\tau)=0$.
We deduce $\lim_{\tau\downarrow 0} v_q(\tau) g_q(\tau)=0$, and therefore  \eqref{2.1.9} is not satisfied.

Next, we consider: (iii) $D=d/\nu$ and $1-\delta D\ge0$; (iv) $D>d/\nu$. 

Using the definitions of $v_q$ and $g_q$, we have
\beqn
v_q(g(\tau)) g_q(g(\tau))= \left[\int_{\tau}^{q(\text{diam}(A))} \left(\log \frac{c}{\rho}\right)^{-\delta D} \rho^{-\nu D+d-1} d\rho\right] \left[\frac{(\tau)^d}{(q(\tau))^D}\right]^{-1}.
\eeqn
Then, computing the limit (for example, applying the L'Hospital's rule), we obtain 
\beqn
\lim_{\tau\downarrow 0} v_q(g(\tau)) g_q(g(\tau)) = (D\nu-d)^{-1}.
\eeqn
Hence, in the case (iii) \eqref{2.1.9} is not satisfied, while in the case (iv) it is.
\medskip

3.\  From \eqref{Lambertequiv} we deduce that the function $\bar g_q$ given in \eqref{2.1.37} satisfies
$\bar g_q(\tau)\asymp \bar g_q(\tau/2)$. Moreover, since $\bar v_q$ and $\bar g_q$ are continuous away from zero, we see that the condition \eqref{2.1.52} is equivalent to 
\beq
\label{2.1.38}
\lim_{\tau\downarrow 0} \bar v_q(\tau) \bar g_q(\tau)\in(0,\infty).
\eeq
Substituting in \eqref{2.1.51} the gauge functions $q_1(\tau)$ and $q_2(\tau)$ by  $\tau^{\nu_1}$ and $\tau^{\nu_2}\left(\log\frac{c}{\tau}\right)^\delta$, respectively, we obtain
\beqn
\bar v_q(\tau) = \nu_1^{-1} \int_\tau^{c_{I,J}} \rho^{-D+\tfrac{d_1}{\nu_1}}\left(q_2^{-1}(\rho)\right)^{d_2-\nu_2}
\left(\log\frac{c}{q_2^{-1}(\rho)}\right)^{1-\delta}\left(\nu_2\log\frac{c}{q_2^{-1}(\rho)} - \delta\right)^{-1}.
\eeqn
Computing the derivative of the reciproque of $\bar g_q$, we see that
\begin{align}
\label{2.1.41}
\frac{d}{d\tau}\left(\left(\bar g_q(\tau)\right)^{-1}\right)
&= \tau^{\tfrac{d_1}{\nu_1}-D-1} \left(q_2^{-1}(\tau)\right)^{d_2-1}
\left[\left(\tfrac{d_1}{\nu_1}-D\right) q_2^{-1}(\tau)\right.\notag\\
&\left.\qquad + d_2 \tau \left(q_2^{-1}(\tau)\right)^{1-\nu_2}
\left(\log\frac{c}{q_2^{-1}(\tau)}\right)^{1-\delta}\left(\nu_2\log\frac{c}{q_2^{-1}(\tau)} - \delta\right)^{-1}\right].
\end{align}
Apply the L'Hospital's rule to obtain
\begin{align*}
\lim_{\tau\downarrow 0}\left[\bar v_q(\tau) \bar g_q(\tau)\right]^{-1} &= \lim_{\tau\downarrow 0}
\frac{\frac{d}{d\tau}\left(\left(\bar g_q(\tau)\right)^{-1}\right)}{\frac{d \bar v_q}{d\tau}(\tau)}
=\lim_{\tau\downarrow 0} \left(L_1(\tau) + L_2(\tau)\right),
\end{align*}
where using \eqref{2.1.41}, we have
\begin{align*}
L_1(\tau)&= \frac{\tau^{\tfrac{d_1}{\nu_1}-D-1} \left(q_2^{-1}(\tau)\right)^{d_2}
\left(\tfrac{d_1}{\nu_1}-D\right)}{\frac{d \bar v_q}{d\tau}(\tau)},\\
L_2(\tau) &= \frac{d_2\tau^{\tfrac{d_1}{\nu_1}-D}\left(q_2^{-1}(\tau)\right)^{d_2-\nu_2}\left(\log\frac{c}{q_2^{-1}(\tau)}\right)^{1-\delta}\left(\nu_2\log\frac{c}{q_2^{-1}(\tau)} - \delta\right)^{-1}}{\frac{d \bar v_q}{d\tau}(\tau)}.
\end{align*}
Since 
$
\frac{1}{2}\nu_2 \log\frac{c}{q_2^{-1}(\tau)} \le \nu_2 \log\frac{c}{q_2^{-1}(\tau)} -\delta \le \nu_2 \log\frac{c}{q_2^{-1}(\tau)}, {\text{as}}\ \tau\downarrow 0, 
$
we find:
\beqn
\lim_{\tau\downarrow 0}L_1(\tau) = D\nu_1\nu_2-d_1\nu_2,\quad \lim_{\tau\downarrow 0}L_2(\tau)= -d_2\nu_1.
\eeqn
Consequently,
\beqn
\lim_{\tau\downarrow 0}\bar v_q(\tau) \bar g_q(\tau) = (D\nu_1\nu_2 - (d_1\nu_2+d_2\nu_1))^{-1}.
\eeqn
This implies \eqref{2.1.38}.

The proof of the lemma is complete.
\qed

\end{document}